\documentclass[11pt]{article}

\usepackage{amssymb,amsmath,amsfonts,epsf,amsmath}
\usepackage{enumerate}

\usepackage[english]{babel}

\usepackage{tikz}

\usetikzlibrary{decorations.pathreplacing,automata,calc,positioning}
\usetikzlibrary{arrows}

\newtheorem{thm}{Theorem}[section]

\newtheorem{cor}[thm]{Corollary}
\newtheorem{lemma}[thm]{Lemma}
\newtheorem{claim}{Claim}
\newtheorem{subclaim}{Claim}[claim]

\newcommand{\efface}[1]{}
\newcommand{\proof}{\noindent{\bf Proof.\ }}

\newcommand{\modo}{{\rm mod \,}}
\newcommand{\even}{{\rm even}}
\newcommand{\odd}{{\rm odd}}
\newcommand{\smallqed}{{\tiny ($\Box$)}}
\newcommand{\QED}{$\Box$}

\newcommand{\cp}{\mathbin{\Box}}
\newcommand{\2}{ \vspace{0.2cm} }
\newcommand{\1}{ \vspace{0.1cm} }

\let\oldenumerate\enumerate
\renewcommand{\enumerate}{
  \oldenumerate
  \setlength{\itemsep}{0pt}
  \setlength{\parskip}{0pt}
  \setlength{\parsep}{0pt}
}

\begin{document}

\title{$3$-Neighbor bootstrap percolation on grids}

\author{$^{1,2}$Jaka Hed\v{z}et \, and \, $^3$Michael A. Henning\thanks{Research supported in part by the University of Johannesburg and the South African National Research Foundation}
\\ \\
Faculty of Natural Sciences and Mathematics \\
University of Maribor, Slovenia \\
$^{2}$Institute of Mathematics, Physics and Mechanics\\
University of Ljubljana, Slovenia \\
\small \tt Email: jaka.hedzet@imfm.com \\ \\
$^{3}$Department of Mathematics and Applied Mathematics \\
University of Johannesburg \\
Auckland Park, 2006 South Africa\\
\small \tt Email: mahenning@uj.ac.za }

\date{}
\maketitle


\begin{abstract}
Given a graph $G$ and assuming that some vertices of $G$ are infected, the $r$-neighbor bootstrap percolation rule makes an uninfected vertex $v$ infected if $v$ has at least $r$ infected neighbors. The $r$-percolation number, $m(G, r)$, of $G$ is the minimum cardinality of a set of initially infected vertices in $G$ such that after continuously performing the $r$-neighbor bootstrap percolation rule each vertex of $G$ eventually becomes infected. In this paper, we consider the $3$-bootstrap percolation number of grids with fixed widths. If $G$ is the cartesian product $P_3 \cp P_m$ of two paths of orders~$3$ and $m$, we prove that $m(G,3)=\frac{3}{2}(m+1)-1$, when $m$ is odd, and $m(G,3)=\frac{3}{2}m +1$, when $m$ is even. Moreover if $G$ is the cartesian product $P_5 \cp P_m$, we prove that $m(G,3)=2m+2$, when $m$ is odd, and $m(G,3)=2m+3$, when $m$ is even. If $G$ is the cartesian product $P_4 \cp P_m$, we prove that $m(G,3)$ takes on one of two possible values, namely $m(G,3) = \lfloor \frac{5(m+1)}{3} \rfloor + 1$ or $m(G,3) = \lfloor \frac{5(m+1)}{3} \rfloor + 2$.
\end{abstract}


\newpage

\section{Introduction}

For notation and graph theory terminology, we in general follow~\cite{HaHeHe-23,HeYe-13}. Specifically, let $G$ be a graph with vertex set $V(G)$ and edge set $E(G)$, and of order~$n(G) = |V(G)|$ and size $m(G) = |E(G)|$. A \emph{neighbor} of a vertex $v$ in $G$ is a vertex $u$ that is adjacent to $v$, that is, $uv \in E(G)$. The \emph{open neighborhood} $N_G(v)$ of a vertex $v$ in $G$ is the set of neighbors of $v$, while the \emph{closed neighborhood} of $v$ is the set $N_G[v] = \{v\} \cup N_G(v)$. For a set $S \subseteq V(G)$, its \emph{open neighborhood} is the set $N_G(S) = \cup_{v \in S} N_G(v)$, and its \emph{closed neighborhood} is the set $N_G[S] = N_G(S) \cup S$.

We denote the \emph{degree} of a vertex $v$ in $G$ by $\deg_G(v)$, or simply by $\deg(v)$ if the graph $G$ is clear from the context, and so $\deg(v) = |N_G(v)|$. If $X \subseteq V(G)$ and $v \in V(G)$, then $\deg_X(v)$ is the number of neighbors of the vertex $v$ in $G$ that belong to the set $X$, that is, $\deg_X(v) = |N_G(v) \cap X|$. In particular, if $X = V(G)$, then $\deg_X(v) = \deg_G(v)$.

We denote a cycle and a path on $n$ vertices by $C_n$ and $P_n$, respectively. For a nonempty set of vertices $S \subseteq V(G)$, the subgraph induced by $S$ is denoted by $G[S]$. Thus, $G[S]$ is the graph having vertex set $S$ and whose edge set consists of all those edges of $G$ incident with two vertices in $S$. Moreover, we denote the graph obtained from $G$ by deleting all vertices in the set $S$ by $G - S$, that is, $G - S = G[V(G) \setminus S]$. A subgraph $H$ of $G$ is an \emph{induced subgraph of $G$} if $H = G[S]$ for some subset $S$ of $V(G)$.

For any integer $r \ge 2$, the $r$-neighbor bootstrap percolation process is an update rule for the states of vertices in a given graph $G$. At any given time a vertex can either be \emph{infected} or \emph{uninfected}. From an initial set of infected vertices, the following updates occur simultaneously and in discrete intervals: any uninfected vertex with at least $r$ infected neighbors becomes infected, while infected vertices never change their state.

More formally, let $A_0 \subseteq V(G)$ be an initial set of infected vertices and for every $t \ge 1$ define
\[
A_t = A_{t-1} \cup \lbrace v\in V(G) \colon |N_G(v) \cap A_{t-1}| \ge r \rbrace.
\]

The set $A_t \setminus A_{t-1}$ is referred to as \emph{the set of vertices infected at time $t$}. A vertex~$v$ is infected before vertex~$u$ if $v \in A_t$ and $u \notin A_t$ for some $t \ge 0$. We say that the set $A_0$ is an \emph{$r$-percolating set}, or simply \emph{$r$-percolates}, in the graph $G$ if
\[
\bigcup\limits_{t=0}^{\infty} A_t= V(G).
\]

A natural extremal problem is finding a smallest $r$-percolating set $A_0$ in a given graph $G$. For a given graph $G$ and integer $r \ge 2$, the \emph{$r$-percolation number} of $G$, denoted $m(G,r)$, is the minimum cardinality of an $r$-percolating set in $G$, that is,
\[
m(G,r)=\min \left\lbrace |A_0| \colon \, A_0 \subseteq V(G), \: \mbox{$A_0$ is an $r$-percolating set in $G$} \right\rbrace.
\]

A \emph{minimum $r$-percolating set} in $G$ is a $r$-percolating set $S$ of $G$ satisfying $m(G,r) = |S|$. Bootstrap percolation is very well studied in graphs, see, for example,~\cite{bal-2006,bal-2012,bal-1998,ben-2021,ben-2014+,bid-2021,bol-2006,
bre-2020,cen-2010,cen-2013,cha-1979,coe-2019,coe-2015,dai-2020,
dnr-2022,gri-2021,gun-2020,mar-2018,mor-2009,prz-2020,prz-2012,romer}.

The \emph{Cartesian product} $G \cp H$ of two graphs $G$ and $H$ is the graph whose vertex set is $V(G) \times V(H)$, and where two vertices $(g_1,h_1)$ and $(g_2,h_2)$ are adjacent in $G \cp H$ if either $g_1=g_2$ and $h_1h_2 \in E(H)$, or $h_1=h_2$ and $g_1g_2 \in E(G)$. For a vertex $g \in V(G)$, the subgraph of $G \cp H$ induced by the set $\{ \, (g,h) \mid h \in V(H) \,\}$ is called a \emph{$H$-fiber} and is denoted by $^g\mspace{1mu}H$. Similarly, for $h \in V(H)$, the \emph{$G$-fiber}, $G^h$, is the subgraph induced by $\{ \, (g,h) \mid g \in V(G) \,\}$. We note that all $G$-fibers are isomorphic to $G$ and all $H$-fibers are isomorphic to~$H$. A \emph{fiber} in $G \cp H$ is a \emph{$G$-fiber} or an \emph{$H$-fiber}.

If $G = P_n \cp  P_m$ is the Cartesian product of two paths $P_n$ and $P_m$ for some $n,m \ge 2$, then a vertex $v \in V(G)$ is called a \emph{boundary vertex} or a \emph{vertex on the boundary of $G$} if $\deg_G(v) \le 3$.

\section{Main results}
\label{S:main}

Our aim in this paper is to study $3$-neighbor bootstrap percolation on grids. We determine closed formulas for the $3$-percolation number of a $3 \times m$ grid for all $m \ge 3$ and the $3$-percolation number of a $5 \times m$ grid for all $m \ge 5$. Moreover, we show that the $3$-percolation number of a $4 \times m$ grid for all $m \ge 4$ takes on one of two possible values. We shall prove the following results.

\begin{thm}
\label{thm:3xm}
For $m \ge 3$, if $G = P_3 \cp P_m$, then
\[
m(G,3)=
\begin{cases}
        \frac{3}{2}(m+1)  - 1; &   m \text{ odd} \1\\
        \frac{3}{2}m + 1 ; &   m \text{ even}
    \end{cases}
\]
\end{thm}

\begin{thm}
\label{thm:5xm}
For $m \ge 5$, if $G = P_5 \cp P_m$, then
\[
m(G,3)=
\begin{cases}
        2m+2; &   m \text{ odd}\\
        2m+3; &   m \text{ even}
    \end{cases}
\]
\end{thm}

\begin{thm}
\label{thm:4xm}
For $m \ge 4$, if $G = P_4 \cp P_m$, then
\[
m(G,3) = \left\lfloor \frac{5(m+1)}{3} \right\rfloor + \Phi_m(G)
\]
where $\Phi_m(G) \in \{1,2\}$. Moreover, $\Phi_m(G) = 1$ if $m \in \{5,7,11\}$.
\end{thm}

\section{Preliminary results}
\label{S:prelim}

In this section, we present some preliminary lemmas that we will need when proving our main results in Section~\ref{S:main}. We remark that Lemma~\ref{lem:subgraph condition} is already known in the literature, but for completeness we provide short proofs of the elementary results we present in this section since we use them frequently when proving our main results.

\begin{lemma}
\label{lem:subgraph condition}
For $r \ge 2$ if $H$ is a subgraph of a graph $G$ such that every vertex in $H$ has strictly less than~$r$ neighbors in $G$ that belong to $V(G) \setminus V(H)$, then every $r$-percolating set of $G$ contains at least one vertex of $H$.
\end{lemma}
\proof For $r \ge 2$ let $H$ be an induced subgraph of a graph $G$ and let $X = V(G) \setminus V(H)$. Suppose that every vertex in $H$ has strictly less than~$r$ neighbors in the graph $G$ that belong to the set $X$, that is, $\deg_X(v) < r$ for every vertex $v \in V(H)$. In this case, even if every vertex in $X$ is infected, no vertex in $H$ becomes infected since every vertex in $H$ has strictly less than $r$ infected neighbors. Therefore, every $r$-percolating set of $G$ contains at least one vertex of $H$.~\QED

\medskip
We call the subgraph $H$ in the statement of Lemma~\ref{lem:subgraph condition} an \emph{$r$-forbidden subgraph of $G$}. We describe next some structural properties of $3$-forbidden subgraphs in grids. As an immediate consequence of Lemma~\ref{lem:subgraph condition}, we infer the following $3$-forbidden subgraphs in grids.

\begin{cor}
\label{cor:subgraph condition_1a}
Let $G = P_n \cp P_m$ for some $m,n \in \mathbb{N}$ and let $S$ be a minimum $3$-percolating set of $G$. If $H$ is a subgraph of $G$ satisfying (a) or (b), then $H$ is a $3$-forbidden subgraph of $G$. \\[-26pt]
\begin{enumerate}
\item[{\rm (a)}] $H$ is a path joining two boundary vertices in $G$;
\item[{\rm (b)}] $H$ is a cycle in $G$.
\end{enumerate}
\end{cor}

We note that if $G = P_n \cp P_m$, then two adjacent boundary vertices in $G$ form a path joining two boundary vertices in $G$. Moreover, every $P_n$-fiber and $P_m$-fiber in $G$ is a path joining two boundary vertices in $G$. We also note that every $4$-cycle in $G$ is a $3$-forbidden subgraph. Hence as special cases of Corollary~\ref{cor:subgraph condition_1a}, we have the following $3$-forbidden subgraphs in a grid.

\begin{cor}
\label{cor:subgraph condition_1}
Let $G = P_n \cp P_m$ for some $m,n \in \mathbb{N}$ and let $S$ be a minimum $3$-percolating set of $G$. If $H$ is an induced subgraph of $G$ satisfying (a), (b) or (c), then $H$ is a $3$-forbidden subgraph of $G$. \\[-26pt]
\begin{enumerate}
\item[{\rm (a)}] $H = P_2$, where the two vertices in $H$ are adjacent boundary vertices in $G$;
\item[{\rm (b)}] $H$ is a fiber in $G$;
\item[{\rm (c)}] $H = C_4$;
\end{enumerate}
\end{cor}

\begin{lemma}
\label{lem:three boundary vertices}
If $G= P_n \cp  P_m$ for some $m,n \in \mathbb{N}$, then there exists a minimum $3$-percolating set of $G$ that does not contain three consecutive boundary vertices of $G$.
\end{lemma}
\proof Let $G= P_n \cp  P_m$ and let $u,v,z$ be three consecutive boundary vertices of $G$ where $v$ is adjacent to both $u$ and $z$. Let $S$ be a minimum $3$-percolating set of $G$ that contains as few vertices from the set $\{u,v,z\}$ as possible. Suppose, to the contrary, that $\{u,v,z\} \subseteq S$. Let $x$ be the third neighbor of $v$ different from $u$ and $z$. If $x \in S$, the $S \setminus \{v\}$ is also a percolating set, since $v$ is adjacent to three infected vertices. However this contradicts the minimality of the set $S$. Therefore, $x \notin S$. We now consider the set $S' = (S \setminus \{v\}) \cup \{x\}$. We note that $|S'| = |S|$. The vertex $v$ becomes immediately infected in the $3$-neighbor bootstrap percolation process since it has three infected neighbors in the set $S'$. Since the resulting set of infected vertices contains the $3$-percolating set $S$ of $G$ as subset, we infer that the set $S'$ is a $3$-percolating set of $G$, implying that $S'$ is a minimum $3$-percolating set of $G$. However since the set $S'$ contains fewer vertices that belong to the set $\{u,v,z\}$ than does the set $S$, this contradicts our choice of the set $S$. Hence, $|\{u,v,z\} \subseteq S| \le 2$, that is, there exists a minimum $3$-percolating set of $G$  that does not contain three consecutive boundary vertices of $G$.~\QED

\subsection{$3$-Bootstrap percolation in $3 \times m$ grids}
\label{S:3xm}

In this section we present a proof of Theorem~\ref{thm:3xm}. Recall its statement.

\noindent \textbf{Theorem~\ref{thm:3xm}} \emph{
For $m \ge 3$, if $G = P_3 \cp P_m$, then
\[
m(G,3)=
\begin{cases}
        \frac{3}{2}(m+1)  - 1; &   m \text{ odd} \1 \\
        \frac{3}{2}m + 1 ; &   m \text{ even}
    \end{cases}
\]
}

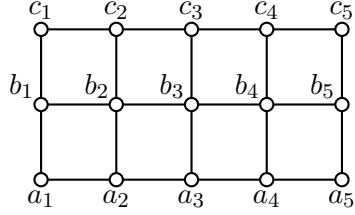
\begin{figure}[htb]
\begin{center}
\begin{tikzpicture}[scale=1,style=thick,x=1cm,y=1cm]
\def\vr{2.5pt} 
\path (0,0) coordinate (A);
\path (0,1) coordinate (B);
\path (0.05,1.25) coordinate (b);
\path (0,2) coordinate (C);
\path (1,0) coordinate (D);
\path (1,1) coordinate (E);
\path (1.05,1.25) coordinate (e);
\path (1,2) coordinate (F);
\path (2,0) coordinate (G);
\path (2,1) coordinate (H);
\path (2.05,1.25) coordinate (h);
\path (2,2) coordinate (I);
\path (3,0) coordinate (J);
\path (3,1) coordinate (K);
\path (3.05,1.25) coordinate (k);
\path (3,2) coordinate (L);
\path (4,0) coordinate (M);
\path (4,1) coordinate (N);
\path (4.05,1.25) coordinate (n);
\path (4,2) coordinate (O);
\draw (A)--(B)--(C);
\draw (D)--(E)--(F);
\draw (G)--(H)--(I);
\draw (J)--(K)--(L);
\draw (M)--(N)--(O);
\draw (A)--(D)--(G)--(J)--(M);
\draw (B)--(E)--(H)--(K)--(N);
\draw (C)--(F)--(I)--(L)--(O);
\draw (A) [fill=white] circle (\vr);
\draw (B) [fill=white] circle (\vr);
\draw (C) [fill=white] circle (\vr);
\draw (D) [fill=white] circle (\vr);
\draw (E) [fill=white] circle (\vr);
\draw (F) [fill=white] circle (\vr);
\draw (G) [fill=white] circle (\vr);
\draw (H) [fill=white] circle (\vr);
\draw (I) [fill=white] circle (\vr);
\draw (J) [fill=white] circle (\vr);
\draw (K) [fill=white] circle (\vr);
\draw (L) [fill=white] circle (\vr);
\draw (M) [fill=white] circle (\vr);
\draw (N) [fill=white] circle (\vr);
\draw (O) [fill=white] circle (\vr);

\draw[anchor = north] (A) node {$a_1$};
\draw[anchor = east] (b) node {$b_1$};
\draw[anchor = south] (C) node {$c_1$};
\draw[anchor = north] (D) node {$a_2$};
\draw[anchor = east] (e) node {$b_2$};
\draw[anchor = south] (F) node {$c_2$};
\draw[anchor = north] (G) node {$a_3$};
\draw[anchor = east] (h) node {$b_3$};
\draw[anchor = south] (I) node {$c_3$};
\draw[anchor = north] (J) node {$a_4$};
\draw[anchor = east] (k) node {$b_4$};
\draw[anchor = south] (L) node {$c_4$};
\draw[anchor = north] (M) node {$a_5$};
\draw[anchor = east] (n) node {$b_5$};
\draw[anchor = south] (O) node {$c_5$};
\end{tikzpicture}
\caption{The graph $G = P_3 \cp  P_5$}
\label{fig:grid1}
\end{center}
\end{figure}

\noindent \textbf{Proof.} For $m \ge 3$, let $G$ be the grid $P_3 \cp P_m$ with
\[
V(G) = \bigcup_{i=1}^{m} \{a_i,b_i,c_i\},
\]
where the path $a_ib_ic_i$ is a $P_3$-fiber in $G$ for $i \in [m]$, and where the paths $a_1a_2 \ldots a_m$, $b_1b_2 \ldots b_m$, and $c_1c_2 \ldots c_m$ are $P_m$-fibers in $G$. For example, when $m = 5$ the grid $G = P_3 \cp P_m$ is illustrated in Figure~\ref{fig:grid1}.

For $i \in [m]$, let $V_i = \{a_i,b_i,c_i\}$ and let
\[
V_{\le i }  = \bigcup_{j=1}^i V_i \hspace*{0.5cm} \mbox{and} \hspace*{0.5cm} V_{\ge i }  = \bigcup_{j=i}^m V_i.
\]

Thus, $V(G) = V_{\le m} = V_{\ge 1}$. Let $A = \{a_1,a_2,\ldots,a_m\}$, $B = \{b_1,b_2,\ldots,b_m\}$, and $C = \{c_1,c_2,\ldots,c_m\}$. Further let
\[
A_\odd = \bigcup_{i=1}^{\lceil \frac{m}{2} \rceil} \{a_{2i-1}\},
\hspace*{0.5cm}
B_\even = \bigcup_{i=1}^{\lfloor \frac{m}{2} \rfloor} \{b_{2i}\},
\hspace*{0.5cm} \mbox{and} \hspace*{0.5cm}
C_\odd = \bigcup_{i=1}^{\lceil \frac{m}{2} \rceil} \{c_{2i-1}\}.
\]

By Lemma~\ref{lem:three boundary vertices}, there exists a minimum $3$-percolating set of $G$ that does not contain three consecutive boundary vertices of $G$. Among all minimum $3$-percolating set of $G$, let $S$ be chosen so that \2 \\
\indent (1) $S$ does not contain three consecutive boundary vertices of $G$, \\
\indent (2) subject to (1), $|S \cap B_\even|$ is a maximum, and \\
\indent (3) subject to (2), $|S \cap (A_\odd \cup C_\odd)|$ is a maximum.

Since each vertex in the set $X = \{a_1, c_1, a_m, c_m\}$ has degree~$2$ in $G$, the set $X$ is necessarily a subset of the $3$-percolating set $S$. Thus since the set $S$ does not contain three consecutive boundary vertices of $G$, we note that
$b_1 \notin S$ and $b_m \notin S$.

Suppose that $m = 3$. In this case, $X = \{a_1, c_1, a_3, c_3\}$. However the set $X$ is not a $3$-percolating set of $G$, implying that $S$ contains at least one additional vertex that does not belong to the set $X$. Since $S \cap \{b_1,a_2,c_2,b_3\} = \emptyset$ by  Lemma~\ref{lem:three boundary vertices}, we infer that $S = X \cup \{b_2\}$, and so $m(G,3) = |S| = 5 = \frac{3}{2}(m+1) - 1$. Hence, we may assume that $m \ge 4$, for otherwise the desired value of $m(G,3)$ holds.

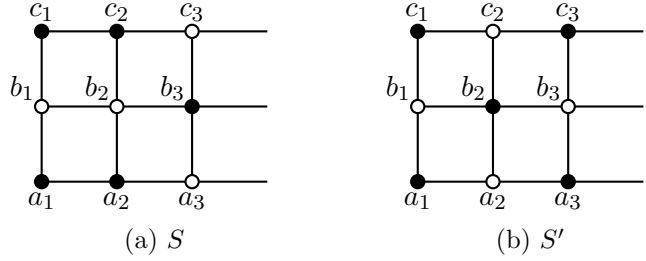
\begin{figure}[htb]
\begin{center}
\begin{tikzpicture}[scale=1,style=thick,x=1cm,y=1cm]
\def\vr{2.5pt} 
\path (0,0) coordinate (A);
\path (0,1) coordinate (B);
\path (0.05,1.25) coordinate (b);
\path (0,2) coordinate (C);
\path (1,0) coordinate (D);
\path (1,1) coordinate (E);
\path (1.05,1.25) coordinate (e);
\path (1,2) coordinate (F);
\path (2,0) coordinate (G);
\path (2,1) coordinate (H);
\path (2.05,1.25) coordinate (h);
\path (2,2) coordinate (I);
\path (3,0) coordinate (J);
\path (3,1) coordinate (K);
\path (3.05,1.25) coordinate (k);
\path (3,2) coordinate (L);
\path (4,0) coordinate (M);
\path (4,1) coordinate (N);
\path (4.05,1.25) coordinate (n);
\path (4,2) coordinate (O);
\path (4,0) coordinate (P);
\path (5,1) coordinate (Q);
\path (5.05,1.25) coordinate (q);
\path (5,2) coordinate (R);
\draw (A)--(B)--(C);
\draw (D)--(E)--(F);
\draw (G)--(H)--(I);
\draw (A)--(D)--(G)--(J);
\draw (B)--(E)--(H)--(K);
\draw (C)--(F)--(I)--(L);

\draw (A) [fill=black] circle (\vr);
\draw (B) [fill=white] circle (\vr);
\draw (C) [fill=black] circle (\vr);
\draw (D) [fill=black] circle (\vr);
\draw (E) [fill=white] circle (\vr);
\draw (F) [fill=black] circle (\vr);
\draw (G) [fill=white] circle (\vr);
\draw (H) [fill=black] circle (\vr);
\draw (I) [fill=white] circle (\vr);
\draw[anchor = north] (A) node {$a_1$};
\draw[anchor = east] (b) node {$b_1$};
\draw[anchor = south] (C) node {$c_1$};
\draw[anchor = north] (D) node {$a_2$};
\draw[anchor = east] (e) node {$b_2$};
\draw[anchor = south] (F) node {$c_2$};
\draw[anchor = north] (G) node {$a_3$};
\draw[anchor = east] (h) node {$b_3$};
\draw[anchor = south] (I) node {$c_3$};
\draw (1.5,-0.8) node {{\small (a) $S$}};
\path (5,0) coordinate (A);
\path (5,1) coordinate (B);
\path (5.05,1.25) coordinate (b);
\path (5,2) coordinate (C);
\path (6,0) coordinate (D);
\path (6,1) coordinate (E);
\path (6.05,1.25) coordinate (e);
\path (6,2) coordinate (F);
\path (7,0) coordinate (G);
\path (7,1) coordinate (H);
\path (7.05,1.25) coordinate (h);
\path (7,2) coordinate (I);
\path (8,0) coordinate (J);
\path (8,1) coordinate (K);
\path (8.05,1.25) coordinate (k);
\path (8,2) coordinate (L);
\draw (A)--(B)--(C);
\draw (D)--(E)--(F);
\draw (G)--(H)--(I);
\draw (A)--(D)--(G)--(J);
\draw (B)--(E)--(H)--(K);
\draw (C)--(F)--(I)--(L);

\draw (A) [fill=black] circle (\vr);
\draw (B) [fill=white] circle (\vr);
\draw (C) [fill=black] circle (\vr);
\draw (D) [fill=white] circle (\vr);
\draw (E) [fill=black] circle (\vr);
\draw (F) [fill=white] circle (\vr);
\draw (G) [fill=black] circle (\vr);
\draw (H) [fill=white] circle (\vr);
\draw (I) [fill=black] circle (\vr);
\draw[anchor = north] (A) node {$a_1$};
\draw[anchor = east] (b) node {$b_1$};
\draw[anchor = south] (C) node {$c_1$};
\draw[anchor = north] (D) node {$a_2$};
\draw[anchor = east] (e) node {$b_2$};
\draw[anchor = south] (F) node {$c_2$};
\draw[anchor = north] (G) node {$a_3$};
\draw[anchor = east] (h) node {$b_3$};
\draw[anchor = south] (I) node {$c_3$};
\draw (6.5,-0.8) node {{\small (b) $S'$}};
\end{tikzpicture}
\caption{The sets $S$ and $S'$ in the proof of Claim~\ref{c:claim1}}
\label{fig:grid2}
\end{center}
\end{figure}

\begin{claim}
\label{c:claim1}
$b_2 \in S$.
\end{claim}
\proof Suppose, to the contrary, that $b_2 \notin S$. Since vertex $b_1$ only gets infected after vertex $b_2$ is infected, the three neighbors $a_2$, $c_2$ and $b_3$ of $b_2$ must all be infected in order to infect $b_2$. Thus, vertex $b_2$ only gets infected after the vertices $a_2$, $c_2$ and $b_3$ are all infected. However if $a_2 \notin S$, then vertex $a_2$ only gets infected after vertex $b_2$ is infected, a contradiction. Hence, $a_2 \in S$. Analogously, $c_2 \in S$. By Lemma~\ref{lem:three boundary vertices}, we infer that $a_3 \notin S$ and $c_3 \notin S$. If $b_3 \notin S$, then it would not be possible to infect~$b_3$ since at most one of its neighbors gets infected. Hence, $b_3 \in S$, and so $S \cap V_{\le 3} = \{a_1,c_1,a_2,c_2,b_3\}$, as illustrated in Figure~\ref{fig:grid2}(a). In this case, we note that the set
\[
S' = (S \setminus \lbrace a_2,c_2,b_3 \rbrace) \cup \lbrace b_2,a_3,c_3 \rbrace
\]
is also a minimum $3$-percolating set of $G$, as illustrated in Figure~\ref{fig:grid2}(b). Thus, $S' \cap V_{\le 3} = \{a_1,c_1,b_2,a_3,c_3\}$ and $S \cap V_{\ge 4} = S' \cap V_{\ge 4}$.

By construction, $|S' \cap B_\even| > |S \cap B_\even|$. Hence if $S'$ satisfies~(1), then we contradict our choice of the set $S$. Therefore, $S'$ does not satisfy~(1), and so $S'$ contains three consecutive boundary vertices. Since the set $S$ does not contain three consecutive boundary vertices, we infer that $\{a_4,a_5\} \subset S$ or $\{c_4,c_5\} \subset S$. If $\{a_4,a_5\} \subset S$ and $\{c_4,c_5\} \subset S$, then the set $S'' = (S' \setminus \{a_4,c_4\}) \cup \{b_4\}$ is a $3$-percolating set of $G$, contradicting the minimality of $S'$. Hence, exactly one of  $\{a_4,a_5\} \subset S$ or $\{c_4,c_5\} \subset S$ holds. By symmetry, we may assume that $\{a_4,a_5\} \subset S$, and so $a_3,a_4,a_5$ are three consecutive boundary vertices that belong to $S'$. The set $S'' = (S' \setminus \{a_4\}) \cup \{b_4\}$ is a minimum $3$-percolating set of $G$ that satisfies~(1). However, $|S'' \cap B_\even| > |S \cap B_\even|$, contradicting our choice of the set $S$. Therefore, $b_2 \in S$.~\smallqed

\medskip
By Claim~\ref{c:claim1}, we have $b_2 \in S$.

\begin{claim}
\label{c:claim2}
$S \cap \{a_2,c_2\} \ne \emptyset$.
\end{claim}
\proof Suppose that at least one of $a_2$ and $c_2$ belongs to the set $S$. By symmetry, we may assume that $a_2 \in S$. In this case, we consider the set $S' = (S \setminus \{a_2\}) \cup \{a_3\}$. Necessarily, $S'$ is a minimum $3$-percolating set of $G$. We note that $|S' \cap B_\even| = |S \cap B_\even|$ and $|S' \cap (A_\odd \cup C_\odd)| > |S \cap (A_\odd \cup C_\odd)|$. If $S'$ satisfies~(1), then we contradict our choice of the set $S$. Hence, $S'$ does not satisfy~(1), implying that $a_3,a_4,a_5$ are three consecutive boundary vertices that belong to $S'$. If $b_4 \in S'$, then $S' \setminus \{a_4\}$ is a  $3$-percolating set of $G$, contradicting the minimality of $S'$. Hence, $b_4 \notin S'$ and the set $S'' = (S' \setminus \{a_4\}) \cup \{b_4\}$ is a minimum $3$-percolating set of $G$ that satisfies~(1) and $|S'' \cap B_\even| > |S \cap B_\even|$, contradicting our choice of the set $S$.~\smallqed

\begin{figure}[htb]
\begin{center}
\begin{tikzpicture}[scale=1,style=thick,x=1cm,y=1cm]
\def\vr{2.5pt} 
\path (0,0) coordinate (A);
\path (0,1) coordinate (B);
\path (0.05,1.25) coordinate (b);
\path (0,2) coordinate (C);
\path (1,0) coordinate (D);
\path (1,1) coordinate (E);
\path (1.05,1.25) coordinate (e);
\path (1,2) coordinate (F);
\path (2,0) coordinate (G);
\path (2,1) coordinate (H);
\path (2.05,1.25) coordinate (h);
\path (2,2) coordinate (I);
\path (3,0) coordinate (J);
\path (3,1) coordinate (K);
\path (3.05,1.25) coordinate (k);
\path (3,2) coordinate (L);
\draw (A)--(B)--(C);
\draw (D)--(E)--(F);
\draw (G)--(H)--(I);
\draw (A)--(D)--(G)--(J);
\draw (B)--(E)--(H)--(K);
\draw (C)--(F)--(I)--(L);

\draw (A) [fill=black] circle (\vr);
\draw (B) [fill=white] circle (\vr);
\draw (C) [fill=black] circle (\vr);
\draw (D) [fill=white] circle (\vr);
\draw (E) [fill=black] circle (\vr);
\draw (F) [fill=white] circle (\vr);
\draw (G) [fill=black] circle (\vr);
\draw (H) [fill=white] circle (\vr);
\draw (I) [fill=black] circle (\vr);
\draw[anchor = north] (A) node {$a_1$};
\draw[anchor = east] (b) node {$b_1$};
\draw[anchor = south] (C) node {$c_1$};
\draw[anchor = north] (D) node {$a_2$};
\draw[anchor = east] (e) node {$b_2$};
\draw[anchor = south] (F) node {$c_2$};
\draw[anchor = north] (G) node {$a_3$};
\draw[anchor = east] (h) node {$b_3$};
\draw[anchor = south] (I) node {$c_3$};
\end{tikzpicture}
\caption{The set $S \cap V_{\le 3}$}
\label{fig:grid3}
\end{center}
\end{figure}
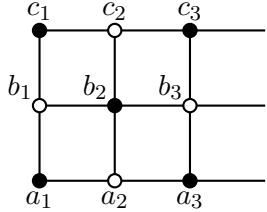

\medskip
By Claim~\ref{c:claim2}, neither $a_2$ nor $c_2$ belongs to the set $S$. Hence, $a_2$ and $c_2$ only get infected after $a_3$ and $c_3$, respectively, are infected. Since $a_3$ and $c_3$ are boundary vertices, this implies that $a_3 \in S$ and $c_3 \in S$. If $b_3 \in S$, then $S \setminus \{b_3\}$ is a $3$-percolating set of $G$, contradicting the minimality of $S$. Hence, $b_3 \notin S$. Thus, $S \cap V_{\le 3} = \{a_1,c_1,b_2,a_3,c_3\}$, as illustrated in Figure~\ref{fig:grid3}.

We note that the set $S$ infects vertices $b_1,a_2,c_2$ and $b_3$. If $m=4$, then by our earlier observations, $\{a_4,c_4\} \subset S$ and $b_4 \notin S$, implying that $S = \{a_1,c_1,b_2,a_3,c_3,a_4,c_4\}$, and so $m(G,3) = |S| = 7 = \frac{3}{2}m + 1$. Hence, we may assume that $m \ge 5$, for otherwise the desired value of $m(G,3)$ holds.

Suppose that $m = 5$. By our earlier observations, $\{a_5,c_5\} \subset S$ and $b_5 \notin S$. In order for $b_5$ to be infected, the vertex $b_4$ must be infected first. If $b_4 \notin S$, then in order for $b_4$ to be infected before $b_5$, both boundary vertices $a_4$ and $c_4$ must belong to the set $S$. But then $(S \setminus \{a_4,c_4\}) \cup \{b_4\}$ is a $3$-percolating set of $G$, contradicting the minimality of $S$. Therefore, $b_4 \in S$, implying that $S = \{a_1,c_1,b_2,a_3,c_3,b_4,a_5,c_5\}$, and so $m(G,3) = |S| = 8 = \frac{3}{2}(m+1) - 1$. Hence, we may assume that $m \ge 6$.

\begin{claim}
\label{c:claim3}
$b_4 \in S$.
\end{claim}
\proof Suppose, to the contrary, that $b_4 \notin S$. For $b_4$ to be infected, it needs two more infected neighbors in addition to the vertex $b_3$, implying that at least one of $a_4$ or $c_4$ must be infected before $b_4$. By symmetry, we may assume that $a_4$ is infected before $b_4$, implying that the boundary vertex $a_4$ belongs to the set $S$ (and to the set $S$). Since $S$ satisfies~(1) and $\{a_3,a_4\} \subset S$, we note that $a_5 \notin S$.

Suppose that $c_4 \notin S$. By Corollary~\ref{cor:subgraph condition_1}(a), the adjacent boundary vertex $c_5$ of $c_4$ therefore belongs to $S$. In order for $c_4$ to be infected, the vertex $b_4$ must be infected first. However in order for $b_4$ to be infected before $c_4$, the vertex $b_5$ must be infected before $b_4$. If $b_5 \notin S$, then both vertices $a_5$ and $b_6$ must be infected before $b_5$, implying that the boundary vertex $a_5$ belongs to the set $S$, a contradiction. Hence, $b_5 \in S$. We now consider the set
\[
S' = (S - \lbrace a_4,b_5 \rbrace ) \cup \lbrace b_4, a_5 \rbrace.
\]

Since $S$ is a $3$-percolating set of $G$, so too is the set $S'$. Thus since $|S'| = |S|$, the set $S'$ is a minimum $3$-percolating set of $G$. We note that $|S' \cap B_\even| > |S \cap B_\even|$, and so if $S'$ satisfies~(1), then we contradict our choice of the set $S$. Hence, $S'$ does not satisfy~(1), implying that $a_5,a_6,a_7$ are three consecutive boundary vertices that belong to $S'$. If $b_6 \in S'$, then $S' \setminus \{a_6\}$ is a  $3$-percolating set of $G$, contradicting the minimality of $S'$. Hence, $b_6 \notin S'$ and the set $S'' = (S' \setminus \{a_6\}) \cup \{b_6\}$ is a minimum $3$-percolating set of $G$ that satisfies~(1) and $|S'' \cap B_\even| > |S \cap B_\even|$, contradicting our choice of the set $S$. Hence, $c_4 \in S$. Since $\{c_3,c_4\} \subset S$ and $S$ satisfies~(1), we note that $c_5 \notin S$.

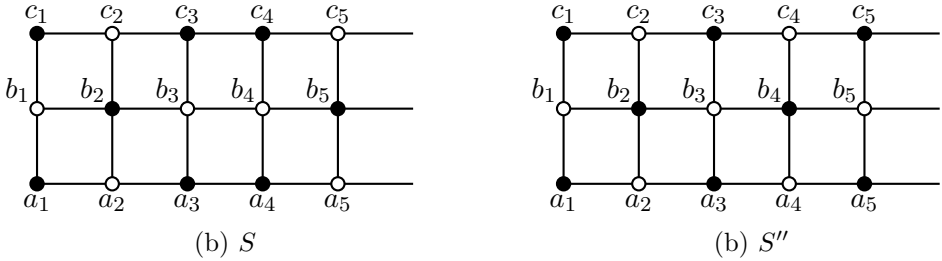
\begin{figure}[htb]
\begin{center}
\begin{tikzpicture}[scale=1,style=thick,x=1cm,y=1cm]
\def\vr{2.5pt} 
\path (0,0) coordinate (A);
\path (0,1) coordinate (B);
\path (0.05,1.25) coordinate (b);
\path (0,2) coordinate (C);
\path (1,0) coordinate (D);
\path (1,1) coordinate (E);
\path (1.05,1.25) coordinate (e);
\path (1,2) coordinate (F);
\path (2,0) coordinate (G);
\path (2,1) coordinate (H);
\path (2.05,1.25) coordinate (h);
\path (2,2) coordinate (I);
\path (3,0) coordinate (J);
\path (3,1) coordinate (K);
\path (3.05,1.25) coordinate (k);
\path (3,2) coordinate (L);
\path (4,0) coordinate (M);
\path (4,1) coordinate (N);
\path (4.05,1.25) coordinate (n);
\path (4,2) coordinate (O);
\path (5,0) coordinate (P);
\path (5,1) coordinate (Q);
\path (5.05,1.25) coordinate (q);
\path (5,2) coordinate (R);
\draw (A)--(B)--(C);
\draw (D)--(E)--(F);
\draw (G)--(H)--(I);
\draw (J)--(K)--(L);
\draw (M)--(N)--(O);
\draw (A)--(D)--(G)--(J)--(M)--(P);
\draw (B)--(E)--(H)--(K)--(N)--(Q);
\draw (C)--(F)--(I)--(L)--(O)--(R);
\draw (A) [fill=black] circle (\vr);
\draw (B) [fill=white] circle (\vr);
\draw (C) [fill=black] circle (\vr);
\draw (D) [fill=white] circle (\vr);
\draw (E) [fill=black] circle (\vr);
\draw (F) [fill=white] circle (\vr);
\draw (G) [fill=black] circle (\vr);
\draw (H) [fill=white] circle (\vr);
\draw (I) [fill=black] circle (\vr);
\draw (J) [fill=black] circle (\vr);
\draw (K) [fill=white] circle (\vr);
\draw (L) [fill=black] circle (\vr);
\draw (M) [fill=white] circle (\vr);
\draw (N) [fill=black] circle (\vr);
\draw (O) [fill=white] circle (\vr);
\draw[anchor = north] (A) node {$a_1$};
\draw[anchor = east] (b) node {$b_1$};
\draw[anchor = south] (C) node {$c_1$};
\draw[anchor = north] (D) node {$a_2$};
\draw[anchor = east] (e) node {$b_2$};
\draw[anchor = south] (F) node {$c_2$};
\draw[anchor = north] (G) node {$a_3$};
\draw[anchor = east] (h) node {$b_3$};
\draw[anchor = south] (I) node {$c_3$};
\draw[anchor = north] (J) node {$a_4$};
\draw[anchor = east] (k) node {$b_4$};
\draw[anchor = south] (L) node {$c_4$};
\draw[anchor = north] (M) node {$a_5$};
\draw[anchor = east] (n) node {$b_5$};
\draw[anchor = south] (O) node {$c_5$};
\draw (2.5,-0.8) node {{\small (b) $S$}};
\path (7,0) coordinate (A);
\path (7,1) coordinate (B);
\path (7.05,1.25) coordinate (b);
\path (7,2) coordinate (C);
\path (8,0) coordinate (D);
\path (8,1) coordinate (E);
\path (8.05,1.25) coordinate (e);
\path (8,2) coordinate (F);
\path (9,0) coordinate (G);
\path (9,1) coordinate (H);
\path (9.05,1.25) coordinate (h);
\path (9,2) coordinate (I);
\path (10,0) coordinate (J);
\path (10,1) coordinate (K);
\path (10.05,1.25) coordinate (k);
\path (10,2) coordinate (L);
\path (11,0) coordinate (M);
\path (11,1) coordinate (N);
\path (11.05,1.25) coordinate (n);
\path (11,2) coordinate (O);
\path (12,0) coordinate (P);
\path (12,1) coordinate (Q);
\path (12.05,1.25) coordinate (q);
\path (12,2) coordinate (R);
\draw (A)--(B)--(C);
\draw (D)--(E)--(F);
\draw (G)--(H)--(I);
\draw (J)--(K)--(L);
\draw (M)--(N)--(O);
\draw (A)--(D)--(G)--(J)--(M)--(P);
\draw (B)--(E)--(H)--(K)--(N)--(Q);
\draw (C)--(F)--(I)--(L)--(O)--(R);
\draw (A) [fill=black] circle (\vr);
\draw (B) [fill=white] circle (\vr);
\draw (C) [fill=black] circle (\vr);
\draw (D) [fill=white] circle (\vr);
\draw (E) [fill=black] circle (\vr);
\draw (F) [fill=white] circle (\vr);
\draw (G) [fill=black] circle (\vr);
\draw (H) [fill=white] circle (\vr);
\draw (I) [fill=black] circle (\vr);
\draw (J) [fill=white] circle (\vr);
\draw (K) [fill=black] circle (\vr);
\draw (L) [fill=white] circle (\vr);
\draw (M) [fill=black] circle (\vr);
\draw (N) [fill=white] circle (\vr);
\draw (O) [fill=black] circle (\vr);
\draw[anchor = north] (A) node {$a_1$};
\draw[anchor = east] (b) node {$b_1$};
\draw[anchor = south] (C) node {$c_1$};
\draw[anchor = north] (D) node {$a_2$};
\draw[anchor = east] (e) node {$b_2$};
\draw[anchor = south] (F) node {$c_2$};
\draw[anchor = north] (G) node {$a_3$};
\draw[anchor = east] (h) node {$b_3$};
\draw[anchor = south] (I) node {$c_3$};
\draw[anchor = north] (J) node {$a_4$};
\draw[anchor = east] (k) node {$b_4$};
\draw[anchor = south] (L) node {$c_4$};
\draw[anchor = north] (M) node {$a_5$};
\draw[anchor = east] (n) node {$b_5$};
\draw[anchor = south] (O) node {$c_5$};
\draw (9.5,-0.8) node {{\small (b) $S''$}};
\end{tikzpicture}
\caption{The sets $S$ and $S''$ in the proof of Claim~\ref{c:claim3}}
\label{fig:grid4}
\end{center}
\end{figure}

If $b_5 \notin S$, then $b_5$ must be infected before the boundary vertices $a_5$ and $c_5$. However, this would not be possible since then $b_5$ would have at most two infected neighbors at any stage of the percolation process. Thus, $b_5 \in S$, and so $S \cap V_{\le 5} = \{a_1,c_1,b_2,a_3,c_3,a_4,c_4,b_5\}$, as illustrated in Figure~\ref{fig:grid4}(a). In this case, we note that the set
\[
S'' = (S \setminus \lbrace a_4,c_4,b_5 \rbrace) \cup \lbrace b_4,a_5,c_5 \rbrace  \1
\]
is also a minimum $3$-percolating set of $G$, as illustrated in Figure~\ref{fig:grid4}(b). Thus, $S' \cap V_{\le 5} = \{a_1,c_1,b_2,a_3,c_3,b_4,a_5,c_5\}$ and $S \cap V_{\ge 6} = S' \cap V_{\ge 6}$.

By construction, $|S'' \cap B_\even| > |S \cap B_\even|$. Hence if $S''$ satisfies~(1), then we contradict our choice of the set $S$. Therefore, $S''$ does not satisfy~(1), and so $S''$ contains three consecutive boundary vertices. Since the set $S$ does not contain three consecutive boundary vertices, we infer that $\{a_6,a_7\} \subset S$ or $\{c_6,c_7\} \subset S$. If $\{a_6,a_7\} \subset S$ and $\{c_6,c_7\} \subset S$, then the set $(S'' \setminus \{a_6,c_6\}) \cup \{b_6\}$ is a $3$-percolating set of $G$, contradicting the minimality of $S''$. Hence, exactly one of $\{a_6,a_7\} \subset S$ or $\{c_6,c_7\} \subset S$ holds. By symmetry, we may assume that $\{a_6,a_7\} \subset S$, and so $a_5,a_6,a_7$ are three consecutive boundary vertices that belong to $S''$. The set $S^* = (S'' \setminus \{a_6\}) \cup \{b_6\}$ is a minimum $3$-percolating set of $G$ that satisfies~(1). However, $|S^* \cap B_\even| > |S \cap B_\even|$, contradicting our choice of the set $S$. Therefore, $b_4 \in S$.~\smallqed

\medskip
By Claim~\ref{c:claim3}, we have $b_4 \in S$.

\begin{claim}
\label{c:claim4}
$S \cap \{a_4,c_4\} \ne \emptyset$.
\end{claim}
\proof Suppose that at least one of $a_4$ and $c_4$ belongs to the set $S$. By symmetry, we may assume that $a_4 \in S$. In this case, we consider the set $S' = (S \setminus \{a_4\}) \cup \{a_5\}$. Necessarily, $S'$ is a minimum $3$-percolating set of $G$. We note that $|S' \cap B_\even| = |S \cap B_\even|$ and $|S' \cap (A_\odd \cup C_\odd)| > |S \cap (A_\odd \cup C_\odd)|$. If $S'$ satisfies~(1), then we contradict our choice of the set $S$. Hence, $S'$ does not satisfy~(1), implying that $a_5,a_6,a_7$ are three consecutive boundary vertices that belong to $S'$. If $b_6 \in S'$, then $S' \setminus \{a_6\}$ is a  $3$-percolating set of $G$, contradicting the minimality of $S'$. Hence, $b_6 \notin S'$ and the set $S'' = (S' \setminus \{a_6\}) \cup \{b_6\}$ is a minimum $3$-percolating set of $G$ that satisfies~(1) and $|S'' \cap B_\even| > |S \cap B_\even|$, contradicting our choice of the set $S$.~\smallqed

\begin{figure}[htb]
\begin{center}
\begin{tikzpicture}[scale=1,style=thick,x=1cm,y=1cm]
\def\vr{2.5pt} 
\path (0,0) coordinate (A);
\path (0,1) coordinate (B);
\path (0.05,1.25) coordinate (b);
\path (0,2) coordinate (C);
\path (1,0) coordinate (D);
\path (1,1) coordinate (E);
\path (1.05,1.25) coordinate (e);
\path (1,2) coordinate (F);
\path (2,0) coordinate (G);
\path (2,1) coordinate (H);
\path (2.05,1.25) coordinate (h);
\path (2,2) coordinate (I);
\path (3,0) coordinate (J);
\path (3,1) coordinate (K);
\path (3.05,1.25) coordinate (k);
\path (3,2) coordinate (L);
\path (4,0) coordinate (M);
\path (4,1) coordinate (N);
\path (4.05,1.25) coordinate (n);
\path (4,2) coordinate (O);
\path (5,0) coordinate (P);
\path (5,1) coordinate (Q);
\path (5.05,1.25) coordinate (q);
\path (5,2) coordinate (R);
\draw (A)--(B)--(C);
\draw (D)--(E)--(F);
\draw (G)--(H)--(I);
\draw (J)--(K)--(L);
\draw (M)--(N)--(O);
\draw (A)--(D)--(G)--(J)--(M)--(P);
\draw (B)--(E)--(H)--(K)--(N)--(Q);
\draw (C)--(F)--(I)--(L)--(O)--(R);
\draw (A) [fill=black] circle (\vr);
\draw (B) [fill=white] circle (\vr);
\draw (C) [fill=black] circle (\vr);
\draw (D) [fill=white] circle (\vr);
\draw (E) [fill=black] circle (\vr);
\draw (F) [fill=white] circle (\vr);
\draw (G) [fill=black] circle (\vr);
\draw (H) [fill=white] circle (\vr);
\draw (I) [fill=black] circle (\vr);
\draw (J) [fill=white] circle (\vr);
\draw (K) [fill=black] circle (\vr);
\draw (L) [fill=white] circle (\vr);
\draw (M) [fill=black] circle (\vr);
\draw (N) [fill=white] circle (\vr);
\draw (O) [fill=black] circle (\vr);
\draw[anchor = north] (A) node {$a_1$};
\draw[anchor = east] (b) node {$b_1$};
\draw[anchor = south] (C) node {$c_1$};
\draw[anchor = north] (D) node {$a_2$};
\draw[anchor = east] (e) node {$b_2$};
\draw[anchor = south] (F) node {$c_2$};
\draw[anchor = north] (G) node {$a_3$};
\draw[anchor = east] (h) node {$b_3$};
\draw[anchor = south] (I) node {$c_3$};
\draw[anchor = north] (J) node {$a_4$};
\draw[anchor = east] (k) node {$b_4$};
\draw[anchor = south] (L) node {$c_4$};
\draw[anchor = north] (M) node {$a_5$};
\draw[anchor = east] (n) node {$b_5$};
\draw[anchor = south] (O) node {$c_5$};
\end{tikzpicture}
\caption{The set $S \cap V_{\le 5}$}
\label{fig:grid5}
\end{center}
\end{figure}
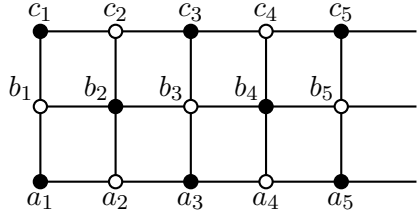

\medskip
By Claim~\ref{c:claim4}, neither $a_4$ nor $c_4$ belongs to the set $S$. Hence, $a_4$ and $c_4$ only get infected after $a_5$ and $c_5$, respectively, are infected. Since $a_5$ and $c_5$ are boundary vertices, this implies that $a_5 \in S$ and $c_5 \in S$. If $b_5 \in S$, then $S \setminus \{b_5\}$ is a $3$-percolating set of $G$, contradicting the minimality of $S$. Hence, $b_5 \notin S$. Thus, $S \cap V_{\le 5} = \{a_1,c_1,b_2,a_3,c_3,b_4,a_5,c_5\}$, as illustrated in Figure~\ref{fig:grid5}.

By our earlier assumption, $m \ge 6$. Continuing the above process, this pattern concludes naturally if $m$ is odd and yields the set
\[
S = A_\odd \cup B_\even \cup C_\odd,
\]
implying that in this case when $m$ is odd, we have
\[
\begin{array}{lcl}
m(G,3) = |S| & = & |A_\odd| + |B_\even| + |C_\odd| \1 \\
& = & \frac{1}{2}(m+1) + \frac{1}{2}(m-1) + \frac{1}{2}(m+1) \1 \\
& = & \frac{3}{2}(m+1) - 1.
\end{array}
\]

If $m$ is even, then recalling that $\{a_m,c_m\} \subset S$ and $b_m \notin S$, this yields the set
\[
S = (A_\odd \cup \{a_m\}) \cup (B_\even \setminus \{b_m\}) \cup (C_\odd \cup \{C_m\}),
\]
implying that in this case when $m$ is even, we have
\[
\begin{array}{lcl}
m(G,3) = |S| & = & (|A_\odd| + 1) + (|B_\even| - 1) + (|C_\odd| + 1) \1 \\
& = & (\frac{1}{2}m + 1) + (\frac{1}{2}m - 1) + (\frac{1}{2}m + 1)  \1 \\
& = & \frac{3}{2}m + 1.
\end{array}
\]

This completes the proof of Theorem~\ref{thm:3xm}.~\QED

\subsection{$3$-Bootstrap percolation in $5 \times m$ grids}
\label{S:5xm}

We present in this section a proof of Theorem~\ref{thm:5xm}. Recall its statement.

\noindent \textbf{Theorem~\ref{thm:5xm}} \emph{
For $m \ge 5$, if $G = P_5 \cp P_m$, then
\[
m(G,3)=
\begin{cases}
        2m+2; &   m \text{ odd}\\
        2m+3; &   m \text{ even}
    \end{cases}
\]
}

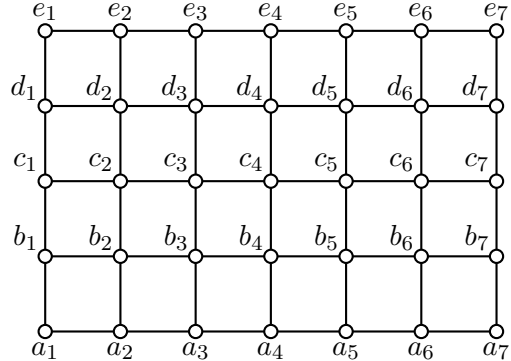
\begin{figure}[htb]
\begin{center}
\begin{tikzpicture}[scale=1,style=thick,x=1cm,y=1cm]
\def\vr{2.5pt} 
\path (0,0) coordinate (A1);
\path (0,1) coordinate (B1);
\path (0.05,1.25) coordinate (b1);
\path (0,2) coordinate (C1);
\path (0.05,2.25) coordinate (c1);
\path (0,3) coordinate (D1);
\path (0.05,3.25) coordinate (d1);
\path (0,4) coordinate (E1);
\path (1,0) coordinate (A2);
\path (1,1) coordinate (B2);
\path (1.05,1.25) coordinate (b2);
\path (1,2) coordinate (C2);
\path (1.05,2.25) coordinate (c2);
\path (1,3) coordinate (D2);
\path (1.05,3.25) coordinate (d2);
\path (1,4) coordinate (E2);
\path (2,0) coordinate (A3);
\path (2,1) coordinate (B3);
\path (2.05,1.25) coordinate (b3);
\path (2,2) coordinate (C3);
\path (2.05,2.25) coordinate (c3);
\path (2,3) coordinate (D3);
\path (2.05,3.25) coordinate (d3);
\path (2,4) coordinate (E3);
\path (3,0) coordinate (A4);
\path (3,1) coordinate (B4);
\path (3.05,1.25) coordinate (b4);
\path (3,2) coordinate (C4);
\path (3.05,2.25) coordinate (c4);
\path (3,3) coordinate (D4);
\path (3.05,3.25) coordinate (d4);
\path (3,4) coordinate (E4);
\path (4,0) coordinate (A5);
\path (4,1) coordinate (B5);
\path (4.05,1.25) coordinate (b5);
\path (4,2) coordinate (C5);
\path (4.05,2.25) coordinate (c5);
\path (4,3) coordinate (D5);
\path (4.05,3.25) coordinate (d5);
\path (4,4) coordinate (E5);
\path (5,0) coordinate (A6);
\path (5,1) coordinate (B6);
\path (5.05,1.25) coordinate (b6);
\path (5,2) coordinate (C6);
\path (5.05,2.25) coordinate (c6);
\path (5,3) coordinate (D6);
\path (5.05,3.25) coordinate (d6);
\path (5,4) coordinate (E6);
\path (6,0) coordinate (A7);
\path (6,1) coordinate (B7);
\path (6.05,1.25) coordinate (b7);
\path (6,2) coordinate (C7);
\path (6.05,2.25) coordinate (c7);
\path (6,3) coordinate (D7);
\path (6.05,3.25) coordinate (d7);
\path (6,4) coordinate (E7);
\draw (A1)--(B1)--(C1)--(D1)--(E1);
\draw (A2)--(B2)--(C2)--(D2)--(E2);
\draw (A3)--(B3)--(C3)--(D3)--(E3);
\draw (A4)--(B4)--(C4)--(D4)--(E4);
\draw (A5)--(B5)--(C5)--(D5)--(E5);
\draw (A6)--(B6)--(C6)--(D6)--(E6);
\draw (A7)--(B7)--(C7)--(D7)--(E7);
\draw (A1)--(A2)--(A3)--(A4)--(A5)--(A6)--(A7);
\draw (B1)--(B2)--(B3)--(B4)--(B5)--(B6)--(B7);
\draw (C1)--(C2)--(C3)--(C4)--(C5)--(C6)--(C7);
\draw (D1)--(D2)--(D3)--(D4)--(D5)--(D6)--(D7);
\draw (E1)--(E2)--(E3)--(E4)--(E5)--(E6)--(E7);
\draw (A1) [fill=white] circle (\vr);
\draw (B1) [fill=white] circle (\vr);
\draw (C1) [fill=white] circle (\vr);
\draw (D1) [fill=white] circle (\vr);
\draw (E1) [fill=white] circle (\vr);
\draw (A2) [fill=white] circle (\vr);
\draw (B2) [fill=white] circle (\vr);
\draw (C2) [fill=white] circle (\vr);
\draw (D2) [fill=white] circle (\vr);
\draw (E2) [fill=white] circle (\vr);
\draw (A3) [fill=white] circle (\vr);
\draw (B3) [fill=white] circle (\vr);
\draw (C3) [fill=white] circle (\vr);
\draw (D3) [fill=white] circle (\vr);
\draw (E3) [fill=white] circle (\vr);
\draw (A4) [fill=white] circle (\vr);
\draw (B4) [fill=white] circle (\vr);
\draw (C4) [fill=white] circle (\vr);
\draw (D4) [fill=white] circle (\vr);
\draw (E4) [fill=white] circle (\vr);
\draw (A5) [fill=white] circle (\vr);
\draw (B5) [fill=white] circle (\vr);
\draw (C5) [fill=white] circle (\vr);
\draw (D5) [fill=white] circle (\vr);
\draw (E5) [fill=white] circle (\vr);
\draw (A6) [fill=white] circle (\vr);
\draw (B6) [fill=white] circle (\vr);
\draw (C6) [fill=white] circle (\vr);
\draw (D6) [fill=white] circle (\vr);
\draw (E6) [fill=white] circle (\vr);
\draw (A7) [fill=white] circle (\vr);
\draw (B7) [fill=white] circle (\vr);
\draw (C7) [fill=white] circle (\vr);
\draw (D7) [fill=white] circle (\vr);
\draw (E7) [fill=white] circle (\vr);
\draw[anchor = north] (A1) node {$a_1$};
\draw[anchor = east] (b1) node {$b_1$};
\draw[anchor = east] (c1) node {$c_1$};
\draw[anchor = east] (d1) node {$d_1$};
\draw[anchor = south] (E1) node {$e_1$};
\draw[anchor = north] (A2) node {$a_2$};
\draw[anchor = east] (b2) node {$b_2$};
\draw[anchor = east] (c2) node {$c_2$};
\draw[anchor = east] (d2) node {$d_2$};
\draw[anchor = south] (E2) node {$e_2$};
\draw[anchor = north] (A3) node {$a_3$};
\draw[anchor = east] (b3) node {$b_3$};
\draw[anchor = east] (c3) node {$c_3$};
\draw[anchor = east] (d3) node {$d_3$};
\draw[anchor = south] (E3) node {$e_3$};
\draw[anchor = north] (A4) node {$a_4$};
\draw[anchor = east] (b4) node {$b_4$};
\draw[anchor = east] (c4) node {$c_4$};
\draw[anchor = east] (d4) node {$d_4$};
\draw[anchor = south] (E4) node {$e_4$};
\draw[anchor = north] (A5) node {$a_5$};
\draw[anchor = east] (b5) node {$b_5$};
\draw[anchor = east] (c5) node {$c_5$};
\draw[anchor = east] (d5) node {$d_5$};
\draw[anchor = south] (E5) node {$e_5$};
\draw[anchor = north] (A6) node {$a_6$};
\draw[anchor = east] (b6) node {$b_6$};
\draw[anchor = east] (c6) node {$c_6$};
\draw[anchor = east] (d6) node {$d_6$};
\draw[anchor = south] (E6) node {$e_6$};
\draw[anchor = north] (A7) node {$a_7$};
\draw[anchor = east] (b7) node {$b_7$};
\draw[anchor = east] (c7) node {$c_7$};
\draw[anchor = east] (d7) node {$d_7$};
\draw[anchor = south] (E7) node {$e_7$};
\end{tikzpicture}
\caption{The graph $G = P_5 \cp  P_7$}
\label{fig:grid6}
\end{center}
\end{figure}

\noindent \textbf{Proof.} For $m \ge 5$, let $G$ be the grid $P_5 \cp P_m$ with
\[
V(G) = \bigcup_{i=1}^{m} \{a_i,b_i,c_i,d_i,e_i\},
\]
where the path $a_ib_ic_id_ie_i$ is a $P_5$-fiber in $G$ for $i \in [m]$, and where $a_1a_2 \ldots a_m$, $b_1b_2 \ldots b_m$, $c_1c_2 \ldots c_m$, $d_1d_2 \ldots d_m$, and $e_1e_2 \ldots e_m$ are $P_m$-fibers in $G$. For example, when $m = 7$ the grid $G = P_5 \cp P_m$ is illustrated in Figure~\ref{fig:grid6}.

For $i \in [m]$, let $V_i = \{a_i,b_i,c_i,d_i,e_i\}$ and let
\[
V_{\le i }  = \bigcup_{j=1}^i V_i \hspace*{0.5cm} \mbox{and} \hspace*{0.5cm} V_{\ge i }  = \bigcup_{j=i}^m V_i.
\]

Thus, $V(G) = V_{\le m} = V_{\ge 1}$. Let $A = \{a_1,a_2,\ldots,a_m\}$, $B = \{b_1,b_2,\ldots,b_m\}$, $C = \{c_1,c_2,\ldots,c_m\}$, $D = \{d_1,d_2,\ldots,d_m\}$, and $E = \{e_1,e_2,\ldots,e_m\}$. In what follows, let $S$ be a minimum $3$-percolating set of $G$ that does not contain three consecutive boundary vertices of $G$. We note that such a set $S$ exists by Lemma~\ref{lem:three boundary vertices}.

\begin{claim}
\label{c:claim5}
The following properties hold. \\ [-26pt]
\begin{enumerate}
\item[{\rm (a)}] $|S \cap V_1| \ge 3$ and $|S \cap (V_1 \cup V_2)| \ge 5$;
\item[{\rm (b)}] $|S \cap V_m| \ge 3$ and $|S \cap (V_{m-1} \cup V_m)| \ge 5$.
\end{enumerate}
\end{claim}
\proof Since the vertices $a_1$ and $e_1$ both have degree~$2$ in $G$, we note that $\{a_1,e_1\} \subset S$. Suppose that $c_1 \notin S$. By Corollary~\ref{cor:subgraph condition_1}(a), this implies that $\{b_1,d_1\} \subset S$. By Corollary~\ref{cor:subgraph condition_1}(b), $|S \cap V_2| \ge 1$. Hence in this case when $c_1 \notin S$, we have $|S \cap V_1| \ge 4$ and $|S \cap (V_1 \cup V_2)| \ge 5$.
Hence we may assume that $c_1 \in S$, for otherwise the desired lower bounds hold. Since $S$ does not contain three consecutive boundary vertices of $G$, we infer that $S \cap V_1 = \{a_1,c_1,e_1\}$. Let $H_1 = G[\{a_2,b_1,b_2\}]$ and let $H_2 = G[\{d_1,d_2,e_2\}]$. Each of $H_1$ and $H_2$ is a path joining two boundary vertices of $G$, and is therefore a $3$-forbidden subgraph of $G$ by Corollary~\ref{cor:subgraph condition_1a}(a). Thus, $S$ contains at least one vertex from each of $H_1$ and $H_2$, implying that $|S \cap \{a_2,b_2\}| \ge 1$ and $|S \cap \{d_2,e_2\}| \ge 1$. Thus, $|S \cap V_1| = 3$ and $|S \cap V_2| \ge 2$, and so $|S \cap (V_1 \cup V_2)| \ge 5$. This proves part~(a). By symmetry, part~(b) holds.~\smallqed

\begin{claim}
\label{c:claim6}
$|S \cap (V_i \cup V_{i+1})| \ge 4$ for all $i$ where $2 \le i \le m-1$.
\end{claim}
\proof  Consider the set $V_i \cup V_{i+1}$ for some $i$ where $2 \le i \le m-1$. We show that $|S \cap (V_i \cup V_{i+1})| \ge 4$. By Corollary~\ref{cor:subgraph condition_1}(a), $|S \cap \{a_{i},a_{i+1}\}| \ge 1$ and $|S \cap\{e_{i},e_{i+1}\}| \ge 1$. By Corollary~\ref{cor:subgraph condition_1}(b), $|S \cap V_i| \ge 1$ and $|S \cap V_{i+1}| \ge 1$.

Suppose firstly that $\{a_{i},e_{i}\} \subset S$. Let $H_1 = G[\{b_i,b_{i+1},c_i,c_{i+1}\}]$ and let $H_2 = G[\{c_i,c_{i+1},d_i,d_{i+1}\}]$. Since $H_1 = C_4$ and $H_2 = C_4$, by Corollary~\ref{cor:subgraph condition_1}(c) both $H_1$ and $H_2$ are $3$-forbidden subgraph of $G$. Let $H_3 = G[V_{i+1}]$. Since $H_3$ is a fibre in $G$, by Corollary~\ref{cor:subgraph condition_1}(b), the fibre $H_3$ is a $3$-forbidden subgraph of $G$. Let $H_4 = G[\{a_{i+1},b_{i+1},b_i,c_i,d_i,d_{i+1},e_{i+1}\}]$. Since $H_4$ is a path joining two boundary vertices, $H_4$ is a $3$-forbidden subgraph of $G$. Hence, the set $S$ must contain at least one vertex from each of the $3$-forbidden subgraphs $H_1$, $H_2$, $H_3$ and $H_4$. This is only possible if $S$ contains at least two vertices in $V_i \cup V_{i+1}$ different from $a_{i}$ and $e_{i}$, implying that $|S \cap (V_i \cup V_{i+1})| \ge 4$.

Hence we may assume that at most one of $a_i$ and $e_i$ belong to the set $S$, for otherwise the desired result holds. By analogous arguments, we may assume that at most one of $a_{i+1}$ and $e_{i+1}$ belong to the set $S$. As observed earlier, at least one $a_i$ and $a_{i+1}$ belongs to the set $S$ and at least one $e_i$ and $e_{i+1}$ belongs to the set $S$. By symmetry, we may therefore assume that $a_{i+1} \in S$ and $e_{i} \in S$, and so $a_i \notin S$ and $e_{i+1} \notin S$.

Let $F_1 = G[\lbrace a_i,b_i,b_{i+1},c_{i+1},d_{i+1},e_{i+1} \rbrace]$, $F_2 = G[\lbrace a_i,b_i,c_i,c_{i+1},d_{i+1},e_{i+1} \rbrace]$ and $F_3 = G[\lbrace a_i,b_i,c_i,d_i,d_{i+1},e_{i+1} \rbrace]$. Each of $F_1$, $F_2$ and $F_3$ is a path joining two boundary vertices of $G$, and is therefore a $3$-forbidden subgraph of $G$ by Corollary~\ref{cor:subgraph condition_1a}(a). Moreover if $F_4 = G[\lbrace b_i,c_i,b_{i+1},c_{i+1} \rbrace]$ and $F_5 = G[\lbrace c_i,d_i,c_{i+1},d_{i+1} \rbrace]$, then $F_4 = C_4$ and $F_5 = C_4$, and so by Corollary~\ref{cor:subgraph condition_1} both $F_4$ and $F_5$ are  $3$-forbidden subgraph of $G$. Hence, the set $S$ must contain at least one vertex from each of the subgraphs $F_1$, $F_2$, $F_3$, $F_4$ and $F_5$. This implies that $S$ contains at least two vertices in $V_i \cup V_{i+1}$ different from $a_{i+1}$ and $e_{i}$, implying once again that $|S \cap (V_i \cup V_{i+1})| \ge 4$.~\smallqed

\newpage
\begin{claim}
\label{c:claim7}
$m(G,3) \ge 2m+2$.
\end{claim}
\proof Suppose firstly that $m$ is even. By Claim~\ref{c:claim5}(a) and Claim~\ref{c:claim6}, we have
\[
\begin{array}{lcl}
m(G,3) = |S| & = & \displaystyle{ \sum_{i=1}^m |S \cap V_i|  }  \1 \\
& = & \displaystyle{ |S \cap V_1| + |S \cap V_m| + \sum_{i=2}^{m-1} |S \cap V_i|   }   \1 \\
& = & \displaystyle{ |S \cap V_1| + |S \cap V_m| + \sum_{i=1}^{\frac{m}{2}-1} |S \cap (V_{2i} \cup V_{2i+1})|  } \1 \\
& \ge & \displaystyle{ 3 + 3 + 4\big(\frac{m}{2}-1 \big)  }  \1 \\
& \ge & \displaystyle{ 2m+2 }.
\end{array}
\]

Suppose secondly that $m$ is odd. By Claim~\ref{c:claim5}(a) and Claim~\ref{c:claim6}, we have
\[
\begin{array}{lcl}
m(G,3) = |S| & = & \displaystyle{ \sum_{i=1}^m |S \cap V_i|  }  \1 \\
& = & \displaystyle{ |S \cap (V_1 \cup V_2)| + |S \cap V_m| + \sum_{i=3}^{m-1} |S \cap V_i|   }   \1 \\
& = & \displaystyle{ |S \cap (V_1 \cup V_2)| + |S \cap V_m| + \sum_{i=1}^{\frac{m-3}{2}} |S \cap (V_{2i+1} \cup V_{2i+2})|  } \1 \\
& \ge & \displaystyle{ 5 + 3 + 4\big( \frac{m-3}{2} \big)  }  \1 \\
& \ge & \displaystyle{ 2m+2 }.
\end{array}
\]
In both cases, we have $m(G,3) \ge 2m+2$.~\smallqed

In order to establish upper bounds on the $3$-percolation number $m(G,3)$ of $G$, let
{\small
\[
A_\odd = \bigcup_{i=1}^{\lceil \frac{m}{2} \rceil} \{a_{2i-1}\},
\hspace*{0.5cm}
B_\even = \bigcup_{i=1}^{\lfloor \frac{m}{2} \rfloor} \{b_{2i}\},
\hspace*{0.5cm}
D_\even = \bigcup_{i=1}^{\lfloor \frac{m}{2} \rfloor} \{d_{2i}\},
\hspace*{0.5cm} \mbox{and} \hspace*{0.5cm}
E_\odd = \bigcup_{i=1}^{\lceil \frac{m}{2} \rceil} \{e_{2i-1}\}.
\]
}

\begin{claim}
\label{c:claim8}
If $m$ is odd, then $m(G,3) = 2m+2$.
\end{claim}
\proof Suppose that $m$ is odd. Let
\[
S_\odd = A_\odd \cup B_\even \cup \{c_1,c_m\} \cup D_\even \cup E_\odd.
\]

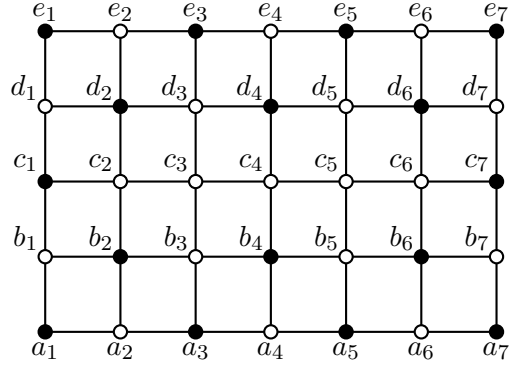
\begin{figure}[htb]
\begin{center}
\begin{tikzpicture}[scale=1,style=thick,x=1cm,y=1cm]
\def\vr{2.5pt} 
\path (0,0) coordinate (A1);
\path (0,1) coordinate (B1);
\path (0.05,1.25) coordinate (b1);
\path (0,2) coordinate (C1);
\path (0.05,2.25) coordinate (c1);
\path (0,3) coordinate (D1);
\path (0.05,3.25) coordinate (d1);
\path (0,4) coordinate (E1);
\path (1,0) coordinate (A2);
\path (1,1) coordinate (B2);
\path (1.05,1.25) coordinate (b2);
\path (1,2) coordinate (C2);
\path (1.05,2.25) coordinate (c2);
\path (1,3) coordinate (D2);
\path (1.05,3.25) coordinate (d2);
\path (1,4) coordinate (E2);
\path (2,0) coordinate (A3);
\path (2,1) coordinate (B3);
\path (2.05,1.25) coordinate (b3);
\path (2,2) coordinate (C3);
\path (2.05,2.25) coordinate (c3);
\path (2,3) coordinate (D3);
\path (2.05,3.25) coordinate (d3);
\path (2,4) coordinate (E3);
\path (3,0) coordinate (A4);
\path (3,1) coordinate (B4);
\path (3.05,1.25) coordinate (b4);
\path (3,2) coordinate (C4);
\path (3.05,2.25) coordinate (c4);
\path (3,3) coordinate (D4);
\path (3.05,3.25) coordinate (d4);
\path (3,4) coordinate (E4);
\path (4,0) coordinate (A5);
\path (4,1) coordinate (B5);
\path (4.05,1.25) coordinate (b5);
\path (4,2) coordinate (C5);
\path (4.05,2.25) coordinate (c5);
\path (4,3) coordinate (D5);
\path (4.05,3.25) coordinate (d5);
\path (4,4) coordinate (E5);
\path (5,0) coordinate (A6);
\path (5,1) coordinate (B6);
\path (5.05,1.25) coordinate (b6);
\path (5,2) coordinate (C6);
\path (5.05,2.25) coordinate (c6);
\path (5,3) coordinate (D6);
\path (5.05,3.25) coordinate (d6);
\path (5,4) coordinate (E6);
\path (6,0) coordinate (A7);
\path (6,1) coordinate (B7);
\path (6.05,1.25) coordinate (b7);
\path (6,2) coordinate (C7);
\path (6.05,2.25) coordinate (c7);
\path (6,3) coordinate (D7);
\path (6.05,3.25) coordinate (d7);
\path (6,4) coordinate (E7);
\draw (A1)--(B1)--(C1)--(D1)--(E1);
\draw (A2)--(B2)--(C2)--(D2)--(E2);
\draw (A3)--(B3)--(C3)--(D3)--(E3);
\draw (A4)--(B4)--(C4)--(D4)--(E4);
\draw (A5)--(B5)--(C5)--(D5)--(E5);
\draw (A6)--(B6)--(C6)--(D6)--(E6);
\draw (A7)--(B7)--(C7)--(D7)--(E7);
\draw (A1)--(A2)--(A3)--(A4)--(A5)--(A6)--(A7);
\draw (B1)--(B2)--(B3)--(B4)--(B5)--(B6)--(B7);
\draw (C1)--(C2)--(C3)--(C4)--(C5)--(C6)--(C7);
\draw (D1)--(D2)--(D3)--(D4)--(D5)--(D6)--(D7);
\draw (E1)--(E2)--(E3)--(E4)--(E5)--(E6)--(E7);
\draw (A1) [fill=black] circle (\vr);
\draw (B1) [fill=white] circle (\vr);
\draw (C1) [fill=black] circle (\vr);
\draw (D1) [fill=white] circle (\vr);
\draw (E1) [fill=black] circle (\vr);
\draw (A2) [fill=white] circle (\vr);
\draw (B2) [fill=black] circle (\vr);
\draw (C2) [fill=white] circle (\vr);
\draw (D2) [fill=black] circle (\vr);
\draw (E2) [fill=white] circle (\vr);
\draw (A3) [fill=black] circle (\vr);
\draw (B3) [fill=white] circle (\vr);
\draw (C3) [fill=white] circle (\vr);
\draw (D3) [fill=white] circle (\vr);
\draw (E3) [fill=black] circle (\vr);
\draw (A4) [fill=white] circle (\vr);
\draw (B4) [fill=black] circle (\vr);
\draw (C4) [fill=white] circle (\vr);
\draw (D4) [fill=black] circle (\vr);
\draw (E4) [fill=white] circle (\vr);
\draw (A5) [fill=black] circle (\vr);
\draw (B5) [fill=white] circle (\vr);
\draw (C5) [fill=white] circle (\vr);
\draw (D5) [fill=white] circle (\vr);
\draw (E5) [fill=black] circle (\vr);
\draw (A6) [fill=white] circle (\vr);
\draw (B6) [fill=black] circle (\vr);
\draw (C6) [fill=white] circle (\vr);
\draw (D6) [fill=black] circle (\vr);
\draw (E6) [fill=white] circle (\vr);
\draw (A7) [fill=black] circle (\vr);
\draw (B7) [fill=white] circle (\vr);
\draw (C7) [fill=black] circle (\vr);
\draw (D7) [fill=white] circle (\vr);
\draw (E7) [fill=black] circle (\vr);
\draw[anchor = north] (A1) node {$a_1$};
\draw[anchor = east] (b1) node {$b_1$};
\draw[anchor = east] (c1) node {$c_1$};
\draw[anchor = east] (d1) node {$d_1$};
\draw[anchor = south] (E1) node {$e_1$};
\draw[anchor = north] (A2) node {$a_2$};
\draw[anchor = east] (b2) node {$b_2$};
\draw[anchor = east] (c2) node {$c_2$};
\draw[anchor = east] (d2) node {$d_2$};
\draw[anchor = south] (E2) node {$e_2$};
\draw[anchor = north] (A3) node {$a_3$};
\draw[anchor = east] (b3) node {$b_3$};
\draw[anchor = east] (c3) node {$c_3$};
\draw[anchor = east] (d3) node {$d_3$};
\draw[anchor = south] (E3) node {$e_3$};
\draw[anchor = north] (A4) node {$a_4$};
\draw[anchor = east] (b4) node {$b_4$};
\draw[anchor = east] (c4) node {$c_4$};
\draw[anchor = east] (d4) node {$d_4$};
\draw[anchor = south] (E4) node {$e_4$};
\draw[anchor = north] (A5) node {$a_5$};
\draw[anchor = east] (b5) node {$b_5$};
\draw[anchor = east] (c5) node {$c_5$};
\draw[anchor = east] (d5) node {$d_5$};
\draw[anchor = south] (E5) node {$e_5$};
\draw[anchor = north] (A6) node {$a_6$};
\draw[anchor = east] (b6) node {$b_6$};
\draw[anchor = east] (c6) node {$c_6$};
\draw[anchor = east] (d6) node {$d_6$};
\draw[anchor = south] (E6) node {$e_6$};
\draw[anchor = north] (A7) node {$a_7$};
\draw[anchor = east] (b7) node {$b_7$};
\draw[anchor = east] (c7) node {$c_7$};
\draw[anchor = east] (d7) node {$d_7$};
\draw[anchor = south] (E7) node {$e_7$};
\end{tikzpicture}
\caption{The set $S_\odd$ in the graph $G = P_5 \cp  P_7$}
\label{fig:grid7}
\end{center}
\end{figure}

For example, when $m = 7$ the set $S_\odd$ is illustrated in Figure~\ref{fig:grid7}.

We show that $S_\odd$ is a $3$-percolating set of $G$. The vertices $b_1,d_1,b_m$ and $d_m$ all have three infected neighbors, and so become infected during the percolation process starting with the initial set $S_\odd$. Hence, all vertices in the first column $V_1$ and in the last column $V_m$ are infected.  The vertices $a_{2i}$ and $e_{2i}$ for $i$ where $1 \le i \le \frac{1}{2}(m-1)$ all have three infected neighbors, and the vertices $b_{2i+1}$ and $d_{2i+1}$ for $i$ where $1 \le i \le \frac{1}{2}(m-3)$ all have three infected neighbors. Therefore, these vertices all become infected during the percolation process. Hence, all vertices in $A \cup B \cup D \cup E$ are infected. Thereafter, all vertices in $C$ become infected by considering the vertices $c_2,c_3,\ldots,c_{m-1}$ sequentially and noting that the vertex $c_i$ becomes infected from its three infected neighbors $c_{i-1}$, $b_i$ and $d_i$ for $i$ where $2 \le i \le m-1$. Hence, all vertices in $V(G)$ become infected, implying that
\[
\begin{array}{lcl}
m(G,3) \le |S_\odd| & = & |A_\odd| + |B_\even| + |\{c_1,c_m\}| + |D_\even| + |E_\odd| \1 \\
& = & \frac{m+1}{2} + \frac{m-1}{2} + 2 + \frac{m-1}{2} + \frac{m+1}{2}  \1 \\
& = & \displaystyle{ 2m+2 }.
\end{array}
\]
Hence, $m(G,3) \le 2m+2$. By Claim~\ref{c:claim7}, $m(G,3) \ge 2m+2$. Consequently, $m(G,3) = 2m+2$.~\smallqed

\begin{claim}
\label{c:claim9}
If $m$ is even, then $m(G,3) \ge 2m + 3$.
\end{claim}
\proof Suppose that $m$ is even. By Claim~\ref{c:claim7}, we know that $m(G,3) \ge 2m + 2$. Suppose, to the contrary, that $m(G,3) = 2m + 2$. Hence we must have equality throughout the inequality chain the first paragraph of the proof of Claim~\ref{c:claim7}, implying that $|S \cap V_1| = |S \cap V_m|= 3$ and $|S \cap (V_{2i} \cup V_{2i+1})| = 4$ for all $i$ where $1 \le i \le \frac{m}{2} - 1$. As shown in the proofs of Claim~\ref{c:claim5} and~\ref{c:claim7}, since $|S \cap V_1| = 3$ we infer that $S \cap V_1 = \{a_1,c_1,e_1\}$ and $|S \cap V_2| \ge 2$. Analogously, since $|S \cap V_m| = 3$ we infer that $S \cap V_m = \{a_m,c_m,e_m\}$ and $|S \cap V_{m-1}| \ge 2$.

\begin{subclaim}
\label{c:claim9.1}
$|S \cap V_i| = 2$ for all $i$ where $2 \le i \le m-1$.
\end{subclaim}
\proof If $|S \cap V_2| \ge 3$, then by our earlier observations we have
{\small
\[
\begin{array}{lcl}
m(G,3) = |S| & = & \displaystyle{ |S \cap V_1| + |S \cap V_2| + \sum_{i=1}^{\frac{m}{2}-2} |S \cap (V_{2i+1} \cup V_{2i+2})| + |S \cap V_{m-1}| + |S \cap V_m|  }   \1 \\
& \ge & \displaystyle{ 3 + 3 + 4\big(\frac{m}{2}-2 \big) + 2 + 3 }  \1 \\
& = & \displaystyle{ 2m+3 },
\end{array}
\]
}

\noindent
a contradiction. Hence, $|S \cap V_2| = 2$. By symmetry, $|S \cap V_{m-1}| = 2$. We show next that $|S \cap V_i| = 2$ for all $i$ where $3 \le i \le m-3$. Let $i$ be the smallest such integer such that $|S \cap V_i| \ne 2$. By Claim~\ref{c:claim6},
$|S \cap (V_{i-1} \cup V_{i})| \ge 4$. By our choice of the integer~$i$, we have $|S \cap V_{i-1}| = 2$, implying that $|S \cap V_i| \ge 3$.

If $i$ is odd, then $|S \cap (V_{i-1} \cup V_{i})| = |S \cap V_{i-1}| + |S \cap V_i| \ge 2 + 3 = 5$, contradicting our earlier observation that $|S \cap (V_{2i} \cup V_{2i+1})| = 4$ for all $i$ where $1 \le i \le \frac{m}{2} - 1$. Hence, $i$ is even. Thus, $4 = |S \cap (V_{i} \cup V_{i+1})| = |S \cap V_{i}| + |S \cap V_{i+1}| \ge 3 + |S \cap V_{i+1}|$, implying that $|S \cap V_{i+1}| = 1$. By Claim~\ref{c:claim6}, $|S \cap (V_{i+1} \cup V_{i+2})| \ge 4$, implying that $|S \cap V_{i+2}| = 3$, which in turn implies that $|S \cap V_{i+3}| = 1$. Continuing in this manner, we have $|S \cap V_{j}| = 3$ for all $j$ even where $i \le j \le m-1$ and $|S \cap V_{j}| = 1$ for all $j$ odd where $i+1 \le j \le m-2$. In particular, $|S \cap V_{m-1}| = 1$, contradicting our earlier observation that $|S \cap V_{m-1}| = 2$.~\smallqed

\medskip
By Claim~\ref{c:claim9.1}, $|S \cap V_i| = 2$ for all $i$ where $2 \le i \le m-1$. By our earlier observation, $S \cap V_1 = \{a_1,c_1,e_1\}$ and $S \cap V_m = \{a_m,c_m,e_m\}$. As shown in the proof of Claim~\ref{c:claim5} we infer that $|S \cap \{a_2,b_2\}| = 1$ and $|S \cap \{d_2,e_2\}| = 1$.  Since $G[\lbrace b_1,b_2,c_2,d_2,d_1 \rbrace]$ is a path joining two boundary vertices of $G$, this subgraph is a $3$-forbidden subgraphs of $G$, implying that $S$ must contain at least one vertex from the set $\{b_2,d_2\}$. By symmetry, we may assume that $b_2 \in S$, and so $a_2 \notin S$. Now either $d_2 \in S$ or $e_2 \in S$. We show firstly that the case $e_2 \in S$ cannot occur.

\begin{subclaim}
\label{c:claim9.2}
If $e_2 \in S$, then we obtain a contradiction.
\end{subclaim}
\proof Suppose that $e_2 \in S$. Thus, $S \cap V_2 = \{b_2,e_2\}$. Since $a_2 \notin S$, the boundary vertex $a_3 \in S$. Let $Q_1 = G[\lbrace c_2,d_2,c_3,d_3 \rbrace]$ and let $Q_2 = G[\lbrace d_1,d_2,d_3,e_3 \rbrace]$. Since $Q_1 = C_4$ and since $Q_2$ is a path joining two boundary vertices of $G$, both $Q_1$ and $Q_2$ are $3$-forbidden subgraphs of $G$, implying that $S$ must contain at least one vertex from each of $Q_1$ and $Q_2$. Since $|S \cap V_3| = 2$ and $a_3 \in S$, we infer that $d_3 \in S$. Thus, $S \cap V_3 = \{a_3,d_3\}$.

Since $e_3 \notin S$, the vertex $e_4 \in S$ by Corollary~\ref{cor:subgraph condition_1}. Let $Q_3 = G[\lbrace b_3,c_3,b_4,c_4 \rbrace]$ and let $Q_2 = G[\lbrace d_1,d_2,c_2,c_3,b_2,b_4,a_4 \rbrace]$. Since $Q_3 = C_4$ and since $Q_4$ is a path joining two boundary vertices of $G$, both $Q_3$ and $Q_4$ are $3$-forbidden subgraphs of $G$, implying that $S$ must contain at least one vertex from each of $Q_3$ and $Q_4$. Since $|S \cap V_4| = 2$ and $e_4 \in S$, we infer that $b_4 \in S$. Thus, $S \cap V_4 = \{b_4,e_4\}$.

Since $a_4 \notin S$, the boundary vertex $a_5 \in S$. Let $Q_5 = G[\lbrace c_4,d_4,c_5,d_5 \rbrace]$ and let $Q_6 = G[\lbrace d_1,d_2,c_2,c_3,c_4,d_4,d_5,e_5 \rbrace]$. Since $Q_5 = C_4$ and since $Q_6$ is a path joining two boundary vertices of $G$, both $Q_5$ and $Q_6$ are $3$-forbidden subgraphs of $G$, implying that $S$ must contain at least one vertex from each of $Q_5$ and $Q_6$. Since $|S \cap V_5| = 2$ and $a_5 \in S$, we infer that $d_5 \in S$. Thus, $S \cap V_5 = \{a_5,d_5\}$.

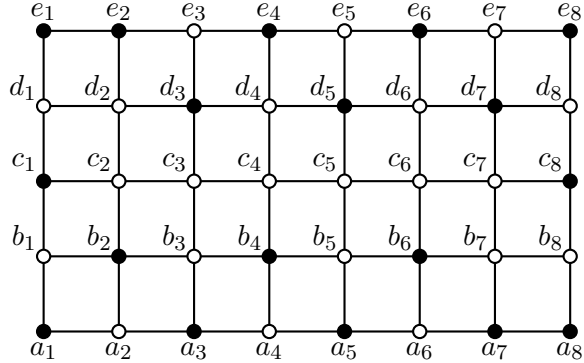
\begin{figure}[htb]
\begin{center}
\begin{tikzpicture}[scale=1,style=thick,x=1cm,y=1cm]
\def\vr{2.5pt} 
\path (0,0) coordinate (A1);
\path (0,1) coordinate (B1);
\path (0.05,1.25) coordinate (b1);
\path (0,2) coordinate (C1);
\path (0.05,2.25) coordinate (c1);
\path (0,3) coordinate (D1);
\path (0.05,3.25) coordinate (d1);
\path (0,4) coordinate (E1);
\path (1,0) coordinate (A2);
\path (1,1) coordinate (B2);
\path (1.05,1.25) coordinate (b2);
\path (1,2) coordinate (C2);
\path (1.05,2.25) coordinate (c2);
\path (1,3) coordinate (D2);
\path (1.05,3.25) coordinate (d2);
\path (1,4) coordinate (E2);
\path (2,0) coordinate (A3);
\path (2,1) coordinate (B3);
\path (2.05,1.25) coordinate (b3);
\path (2,2) coordinate (C3);
\path (2.05,2.25) coordinate (c3);
\path (2,3) coordinate (D3);
\path (2.05,3.25) coordinate (d3);
\path (2,4) coordinate (E3);
\path (3,0) coordinate (A4);
\path (3,1) coordinate (B4);
\path (3.05,1.25) coordinate (b4);
\path (3,2) coordinate (C4);
\path (3.05,2.25) coordinate (c4);
\path (3,3) coordinate (D4);
\path (3.05,3.25) coordinate (d4);
\path (3,4) coordinate (E4);
\path (4,0) coordinate (A5);
\path (4,1) coordinate (B5);
\path (4.05,1.25) coordinate (b5);
\path (4,2) coordinate (C5);
\path (4.05,2.25) coordinate (c5);
\path (4,3) coordinate (D5);
\path (4.05,3.25) coordinate (d5);
\path (4,4) coordinate (E5);
\path (5,0) coordinate (A6);
\path (5,1) coordinate (B6);
\path (5.05,1.25) coordinate (b6);
\path (5,2) coordinate (C6);
\path (5.05,2.25) coordinate (c6);
\path (5,3) coordinate (D6);
\path (5.05,3.25) coordinate (d6);
\path (5,4) coordinate (E6);
\path (6,0) coordinate (A7);
\path (6,1) coordinate (B7);
\path (6.05,1.25) coordinate (b7);
\path (6,2) coordinate (C7);
\path (6.05,2.25) coordinate (c7);
\path (6,3) coordinate (D7);
\path (6.05,3.25) coordinate (d7);
\path (6,4) coordinate (E7);
\path (7,0) coordinate (A8);
\path (7,1) coordinate (B8);
\path (7.05,1.25) coordinate (b8);
\path (7,2) coordinate (C8);
\path (7.05,2.25) coordinate (c8);
\path (7,3) coordinate (D8);
\path (7.05,3.25) coordinate (d8);
\path (7,4) coordinate (E8);
\draw (A1)--(B1)--(C1)--(D1)--(E1);
\draw (A2)--(B2)--(C2)--(D2)--(E2);
\draw (A3)--(B3)--(C3)--(D3)--(E3);
\draw (A4)--(B4)--(C4)--(D4)--(E4);
\draw (A5)--(B5)--(C5)--(D5)--(E5);
\draw (A6)--(B6)--(C6)--(D6)--(E6);
\draw (A7)--(B7)--(C7)--(D7)--(E7);
\draw (A8)--(B8)--(C8)--(D8)--(E8);
\draw (A1)--(A2)--(A3)--(A4)--(A5)--(A6)--(A7)--(A8);
\draw (B1)--(B2)--(B3)--(B4)--(B5)--(B6)--(B7)--(B8);
\draw (C1)--(C2)--(C3)--(C4)--(C5)--(C6)--(C7)--(C8);
\draw (D1)--(D2)--(D3)--(D4)--(D5)--(D6)--(D7)--(D8);
\draw (E1)--(E2)--(E3)--(E4)--(E5)--(E6)--(E7)--(E8);
\draw (A1) [fill=black] circle (\vr);
\draw (B1) [fill=white] circle (\vr);
\draw (C1) [fill=black] circle (\vr);
\draw (D1) [fill=white] circle (\vr);
\draw (E1) [fill=black] circle (\vr);
\draw (A2) [fill=white] circle (\vr);
\draw (B2) [fill=black] circle (\vr);
\draw (C2) [fill=white] circle (\vr);
\draw (D2) [fill=white] circle (\vr);
\draw (E2) [fill=black] circle (\vr);
\draw (A3) [fill=black] circle (\vr);
\draw (B3) [fill=white] circle (\vr);
\draw (C3) [fill=white] circle (\vr);
\draw (D3) [fill=black] circle (\vr);
\draw (E3) [fill=white] circle (\vr);
\draw (A4) [fill=white] circle (\vr);
\draw (B4) [fill=black] circle (\vr);
\draw (C4) [fill=white] circle (\vr);
\draw (D4) [fill=white] circle (\vr);
\draw (E4) [fill=black] circle (\vr);
\draw (A5) [fill=black] circle (\vr);
\draw (B5) [fill=white] circle (\vr);
\draw (C5) [fill=white] circle (\vr);
\draw (D5) [fill=black] circle (\vr);
\draw (E5) [fill=white] circle (\vr);
\draw (A6) [fill=white] circle (\vr);
\draw (B6) [fill=black] circle (\vr);
\draw (C6) [fill=white] circle (\vr);
\draw (D6) [fill=white] circle (\vr);
\draw (E6) [fill=black] circle (\vr);
\draw (A7) [fill=black] circle (\vr);
\draw (B7) [fill=white] circle (\vr);
\draw (C7) [fill=white] circle (\vr);
\draw (D7) [fill=black] circle (\vr);
\draw (E7) [fill=white] circle (\vr);
\draw (A8) [fill=black] circle (\vr);
\draw (B8) [fill=white] circle (\vr);
\draw (C8) [fill=black] circle (\vr);
\draw (D8) [fill=white] circle (\vr);
\draw (E8) [fill=black] circle (\vr);
\draw[anchor = north] (A1) node {$a_1$};
\draw[anchor = east] (b1) node {$b_1$};
\draw[anchor = east] (c1) node {$c_1$};
\draw[anchor = east] (d1) node {$d_1$};
\draw[anchor = south] (E1) node {$e_1$};
\draw[anchor = north] (A2) node {$a_2$};
\draw[anchor = east] (b2) node {$b_2$};
\draw[anchor = east] (c2) node {$c_2$};
\draw[anchor = east] (d2) node {$d_2$};
\draw[anchor = south] (E2) node {$e_2$};
\draw[anchor = north] (A3) node {$a_3$};
\draw[anchor = east] (b3) node {$b_3$};
\draw[anchor = east] (c3) node {$c_3$};
\draw[anchor = east] (d3) node {$d_3$};
\draw[anchor = south] (E3) node {$e_3$};
\draw[anchor = north] (A4) node {$a_4$};
\draw[anchor = east] (b4) node {$b_4$};
\draw[anchor = east] (c4) node {$c_4$};
\draw[anchor = east] (d4) node {$d_4$};
\draw[anchor = south] (E4) node {$e_4$};
\draw[anchor = north] (A5) node {$a_5$};
\draw[anchor = east] (b5) node {$b_5$};
\draw[anchor = east] (c5) node {$c_5$};
\draw[anchor = east] (d5) node {$d_5$};
\draw[anchor = south] (E5) node {$e_5$};
\draw[anchor = north] (A6) node {$a_6$};
\draw[anchor = east] (b6) node {$b_6$};
\draw[anchor = east] (c6) node {$c_6$};
\draw[anchor = east] (d6) node {$d_6$};
\draw[anchor = south] (E6) node {$e_6$};
\draw[anchor = north] (A7) node {$a_7$};
\draw[anchor = east] (b7) node {$b_7$};
\draw[anchor = east] (c7) node {$c_7$};
\draw[anchor = east] (d7) node {$d_7$};
\draw[anchor = south] (E7) node {$e_7$};
\draw[anchor = north] (A8) node {$a_8$};
\draw[anchor = east] (b8) node {$b_8$};
\draw[anchor = east] (c8) node {$c_8$};
\draw[anchor = east] (d8) node {$d_8$};
\draw[anchor = south] (E8) node {$e_8$};
\end{tikzpicture}
\caption{The set $S$ in the graph $G = P_5 \cp  P_8$ in the proof of Claim~\ref{c:claim9.2}}
\label{fig:grid8}
\end{center}
\end{figure}

Continuing in this way, the above pattern repeats itself, that is, $S \cap V_i = \{b_i,e_i\}$ for $i$ even and $2 \le i \le m-2$ and $S \cap V_i = \{a_i,d_i\}$ for $i$ odd and $3 \le i \le m-1$. The set $S$ is now fully determined. For example, when $m = 8$ the set $S$ is illustrated in Figure~\ref{fig:grid8}. However, the subgraph $G[\lbrace \{d_1,d_m\} \cup (V(C) \setminus \{c_1,c_m\} \rbrace]$ is a path joining two boundary vertices of $G$ and is therefore a $3$-forbidden subgraph of $G$. However, this subgraph contains no vertex of $S$, a contradiction.~\smallqed

\medskip
By Claim~\ref{c:claim9.2}, $e_2 \notin S$. By our earlier observations, $|S \cap \{d_2,e_2\}| = 1$, implying that $d_2 \in S$. Thus, $S \cap V_2 = \{b_2,d_2\}$. Since $a_2 \notin S$, this forces $a_3 \in S$, and since $e_2 \notin S$, this forces $e_3 \in S$. Thus, $S \cap V_3 = \{a_3,e_3\}$.

Let $R_1 = G[\lbrace b_3,b_4,c_3,c_4 \rbrace]$ and let $R_2 = G[\lbrace c_3,c_4,d_3,d_4 \rbrace]$. Since $R_1 = C_4$ and  $R_1 = C_4$, both $R_1$ and $R_2$ are $3$-forbidden subgraphs of $G$, implying that $S$ must contain at least one vertex from each of $R_1$ and $R_2$. This implies that at most one of $a_4$ and $e_4$ belong to the set $S$. By symmetry, we may assume that $e_4 \notin S$, implying that $e_5 \in S$. If $a_4 \in S$, then this forces $c_4 \in S$ in order for the set $S$ to contain a vertex from each of $R_1$ and $R_2$. We note that the case $S \cap V_4 = \{b_4,c_4\}$ is symmetric to the case $S \cap V_4 = \{c_4,d_4\}$. Hence by symmetry, there are three possibilities for the set $S \cap V_4$, namely $S \cap V_4 = \{b_4,c_4\}$, $S \cap V_4 = \{b_4,d_4\}$, or $S \cap V_4 = \{a_4,c_4\}$.

We show next that the cases $S \cap V_4 = \{a_4,c_4\}$ and $S \cap V_4 = \{b_4,c_4\}$ cannot occur.

\begin{subclaim}
\label{c:claim9.3}
If $S \cap V_4 = \{a_4,c_4\}$, then we obtain a contradiction.
\end{subclaim}
\proof Suppose that $S \cap V_4 = \{a_4,c_4\}$. Since $\{a_3,a_4\} \subset S$, we know that $a_5 \notin S$. Let $L_1 = G[\lbrace e_4,d_4,d_3,c_3,b_3,b_4,b_5,a_5 \rbrace]$. Since $L_1$ is a path joining two boundary vertices of $G$, the subgraph $L_1$ is a $3$-forbidden subgraphs of $G$, and so $S$ must contain at least one vertex from $L_1$, implying that $b_5 \in S$. Thus, $S \cap V_5 = \{b_5,e_5\}$.

Since $a_5 \notin S$, this forces $a_6 \in S$. Let $L_2 = G[\lbrace e_4,d_4,d_5,d_6,e_6 \rbrace]$ and let $L_3 = G[\lbrace c_5,c_6,d_5,d_6 \rbrace]$. Since $L_2$ is a path joining two boundary vertices of $G$ and since $L_3 = C_4$, the subgraphs $L_2$ and $L_3$ are $3$-forbidden subgraphs of $G$, and so $S$ must contain at least one vertex from each of $L_2$ and $L_3$, implying that $d_6 \in S$. Thus, $S \cap V_6 = \{a_6,d_6\}$.

Since $e_6 \notin S$, this forces $e_7 \in S$. Let $L_4 = G[\lbrace e_4,d_4,d_5,c_5,c_6,b_6,b_7,a_7 \rbrace]$ and let $L_5 = G[\lbrace b_6,b_7,c_6,c_7 \rbrace]$. Since $L_4$ is a path joining two boundary vertices of $G$ and since $L_5 = C_4$, the subgraphs $L_4$ and $L_5$ are $3$-forbidden subgraphs of $G$, and so $S$ must contain at least one vertex from each of $L_4$ and $L_5$, implying that $b_7 \in S$. Thus, $S \cap V_7 = \{b_7,e_7\}$.

Continuing in this way, the above pattern repeats itself, that is, $S \cap V_i = \{b_i,e_i\}$ for $i$ odd and $5 \le i \le m-1$ and $S \cap V_i = \{a_i,d_i\}$ for $i$ even and $6 \le i \le m-2$. The set $S$ is now fully determined. However, the subgraph $G[\{e_4,d_4,d_5,d_{m-1},d_m\} \cup (V(C) \setminus \{c_1,c_2,c_3,c_4,c_m\})]$ is a path joining two boundary vertices of $G$ and is therefore a $3$-forbidden subgraph of $G$. However, this subgraph contains no vertex of $S$, a contradiction.~\smallqed

\begin{subclaim}
\label{c:claim9.4}
If $S \cap V_4 = \{b_4,c_4\}$, then we obtain a contradiction.
\end{subclaim}
\proof Suppose that $S \cap V_4 = \{b_4,c_4\}$. Since $a_4 \notin S$, this forces $a_5 \in S$. Recall that $e_5 \in S$, and so $S \cap V_5 = \{a_5,e_5\}$. Let $Z = \{b_m,b_{m-1},c_{m-1},d_{m-1},d_m\}$.

If $m = 6$, then the set $S$ is fully determined. In this case, the subgraph $G[Z]$ is a path joining two boundary vertices of $G$ and is therefore a $3$-forbidden subgraph of $G$. However, this subgraph contains no vertex of $S$, a contradiction. Hence, $m \ge 8$. Let $T_1 = G[\lbrace e_4,d_4,d_5,d_6,e_6\rbrace]$, $T_2 = G[\lbrace c_5,c_6,d_5,d_6 \rbrace]$, $T_3 = G[\lbrace e_4,d_4,d_5,c_5,b_5,b_6,a_6 \rbrace]$, and $T_4 = G[\lbrace b_5,b_6,c_5,c_6 \rbrace]$. Since $T_1$ and $T_3$ are paths joining two boundary vertices of $G$ and since $T_2 = C_4$ and $T_4 = C_4$, the subgraphs $T_1, T_2, T_3$ and $T_4$ are $3$-forbidden subgraphs of $G$, and so $S$ must contain at least one vertex from each of $T_1, T_2, T_3$ and $T_4$, implying that $S \cap V_6 = \{b_6,d_6\}$. Since $a_6 \notin S$, this forces $a_7 \in S$, and since $e_6 \notin S$, this forces $e_7 \in S$, and so $S \cap V_7 = \{a_6,e_6\}$.

Continuing in this way, the above pattern repeats itself, that is, $S \cap V_i = \{a_i,e_i\}$ for $i$ odd and $5 \le i \le m-1$ and $S \cap V_i = \{b_i,d_i\}$ for $i$ even and $6 \le i \le m-2$. The set $S$ is now fully determined. However, as before the subgraph $G[Z]$ is a path joining two boundary vertices of $G$ and is therefore a $3$-forbidden subgraph of $G$. However, this subgraph contains no vertex of $S$, a contradiction.~\smallqed

\medskip
By Claim~\ref{c:claim9.3}, the case $S \cap V_4 = \{a_4,c_4\}$ cannot occur. By Claim~\ref{c:claim9.4}, the case $S \cap V_4 = \{b_4,c_4\}$ cannot occur. Hence by our earlier assumptions, $S \cap V_4 = \{b_4,d_4\}$. Since $a_4 \notin S$, this forces $a_5 \in S$, and since $e_4 \notin S$, this forces $e_5 \in S$, and so $S \cap V_5 = \{a_5,e_5\}$. Let $Z = \{b_m,b_{m-1},c_{m-1},d_{m-1},d_m\}$.

If $m = 6$, then the set $S$ is fully determined. In this case, the subgraph $G[Z]$ is a path joining two boundary vertices of $G$ and is therefore a $3$-forbidden subgraph of $G$. However, this subgraph contains no vertex of $S$, a contradiction. Hence, $m \ge 8$.

If $S \cap V_6 = \{a_6,c_6\}$, then proceeding analogously as in the proof of Claim~\ref{c:claim9.3} we obtain a contradiction. If $S \cap V_6 = \{b_6,c_6\}$, then proceeding analogously as in the proof of Claim~\ref{c:claim9.4} we obtain a contradiction. Hence, $S \cap V_6 = \{b_6,d_6\}$.

\begin{figure}[htb]
\begin{center}
\begin{tikzpicture}[scale=1,style=thick,x=1cm,y=1cm]
\def\vr{2.5pt} 
\path (0,0) coordinate (A1);
\path (0,1) coordinate (B1);
\path (0.05,1.25) coordinate (b1);
\path (0,2) coordinate (C1);
\path (0.05,2.25) coordinate (c1);
\path (0,3) coordinate (D1);
\path (0.05,3.25) coordinate (d1);
\path (0,4) coordinate (E1);
\path (1,0) coordinate (A2);
\path (1,1) coordinate (B2);
\path (1.05,1.25) coordinate (b2);
\path (1,2) coordinate (C2);
\path (1.05,2.25) coordinate (c2);
\path (1,3) coordinate (D2);
\path (1.05,3.25) coordinate (d2);
\path (1,4) coordinate (E2);
\path (2,0) coordinate (A3);
\path (2,1) coordinate (B3);
\path (2.05,1.25) coordinate (b3);
\path (2,2) coordinate (C3);
\path (2.05,2.25) coordinate (c3);
\path (2,3) coordinate (D3);
\path (2.05,3.25) coordinate (d3);
\path (2,4) coordinate (E3);
\path (3,0) coordinate (A4);
\path (3,1) coordinate (B4);
\path (3.05,1.25) coordinate (b4);
\path (3,2) coordinate (C4);
\path (3.05,2.25) coordinate (c4);
\path (3,3) coordinate (D4);
\path (3.05,3.25) coordinate (d4);
\path (3,4) coordinate (E4);
\path (4,0) coordinate (A5);
\path (4,1) coordinate (B5);
\path (4.05,1.25) coordinate (b5);
\path (4,2) coordinate (C5);
\path (4.05,2.25) coordinate (c5);
\path (4,3) coordinate (D5);
\path (4.05,3.25) coordinate (d5);
\path (4,4) coordinate (E5);
\path (5,0) coordinate (A6);
\path (5,1) coordinate (B6);
\path (5.05,1.25) coordinate (b6);
\path (5,2) coordinate (C6);
\path (5.05,2.25) coordinate (c6);
\path (5,3) coordinate (D6);
\path (5.05,3.25) coordinate (d6);
\path (5,4) coordinate (E6);
\path (6,0) coordinate (A7);
\path (6,1) coordinate (B7);
\path (6.05,1.25) coordinate (b7);
\path (6,2) coordinate (C7);
\path (6.05,2.25) coordinate (c7);
\path (6,3) coordinate (D7);
\path (6.05,3.25) coordinate (d7);
\path (6,4) coordinate (E7);
\path (7,0) coordinate (A8);
\path (7,1) coordinate (B8);
\path (7.05,1.25) coordinate (b8);
\path (7,2) coordinate (C8);
\path (7.05,2.25) coordinate (c8);
\path (7,3) coordinate (D8);
\path (7.05,3.25) coordinate (d8);
\path (7,4) coordinate (E8);
\draw (A1)--(B1)--(C1)--(D1)--(E1);
\draw (A2)--(B2)--(C2)--(D2)--(E2);
\draw (A3)--(B3)--(C3)--(D3)--(E3);
\draw (A4)--(B4)--(C4)--(D4)--(E4);
\draw (A5)--(B5)--(C5)--(D5)--(E5);
\draw (A6)--(B6)--(C6)--(D6)--(E6);
\draw (A7)--(B7)--(C7)--(D7)--(E7);
\draw (A8)--(B8)--(C8)--(D8)--(E8);
\draw (A1)--(A2)--(A3)--(A4)--(A5)--(A6)--(A7)--(A8);
\draw (B1)--(B2)--(B3)--(B4)--(B5)--(B6)--(B7)--(B8);
\draw (C1)--(C2)--(C3)--(C4)--(C5)--(C6)--(C7)--(C8);
\draw (D1)--(D2)--(D3)--(D4)--(D5)--(D6)--(D7)--(D8);
\draw (E1)--(E2)--(E3)--(E4)--(E5)--(E6)--(E7)--(E8);
\draw (A1) [fill=black] circle (\vr);
\draw (B1) [fill=white] circle (\vr);
\draw (C1) [fill=black] circle (\vr);
\draw (D1) [fill=white] circle (\vr);
\draw (E1) [fill=black] circle (\vr);
\draw (A2) [fill=white] circle (\vr);
\draw (B2) [fill=black] circle (\vr);
\draw (C2) [fill=white] circle (\vr);
\draw (D2) [fill=black] circle (\vr);
\draw (E2) [fill=white] circle (\vr);
\draw (A3) [fill=black] circle (\vr);
\draw (B3) [fill=white] circle (\vr);
\draw (C3) [fill=white] circle (\vr);
\draw (D3) [fill=white] circle (\vr);
\draw (E3) [fill=black] circle (\vr);
\draw (A4) [fill=white] circle (\vr);
\draw (B4) [fill=black] circle (\vr);
\draw (C4) [fill=white] circle (\vr);
\draw (D4) [fill=black] circle (\vr);
\draw (E4) [fill=white] circle (\vr);
\draw (A5) [fill=black] circle (\vr);
\draw (B5) [fill=white] circle (\vr);
\draw (C5) [fill=white] circle (\vr);
\draw (D5) [fill=white] circle (\vr);
\draw (E5) [fill=black] circle (\vr);
\draw (A6) [fill=white] circle (\vr);
\draw (B6) [fill=black] circle (\vr);
\draw (C6) [fill=white] circle (\vr);
\draw (D6) [fill=black] circle (\vr);
\draw (E6) [fill=white] circle (\vr);
\draw (A7) [fill=black] circle (\vr);
\draw (B7) [fill=white] circle (\vr);
\draw (C7) [fill=white] circle (\vr);
\draw (D7) [fill=white] circle (\vr);
\draw (E7) [fill=black] circle (\vr);
\draw (A8) [fill=black] circle (\vr);
\draw (B8) [fill=white] circle (\vr);
\draw (C8) [fill=black] circle (\vr);
\draw (D8) [fill=white] circle (\vr);
\draw (E8) [fill=black] circle (\vr);
\draw[anchor = north] (A1) node {$a_1$};
\draw[anchor = east] (b1) node {$b_1$};
\draw[anchor = east] (c1) node {$c_1$};
\draw[anchor = east] (d1) node {$d_1$};
\draw[anchor = south] (E1) node {$e_1$};
\draw[anchor = north] (A2) node {$a_2$};
\draw[anchor = east] (b2) node {$b_2$};
\draw[anchor = east] (c2) node {$c_2$};
\draw[anchor = east] (d2) node {$d_2$};
\draw[anchor = south] (E2) node {$e_2$};
\draw[anchor = north] (A3) node {$a_3$};
\draw[anchor = east] (b3) node {$b_3$};
\draw[anchor = east] (c3) node {$c_3$};
\draw[anchor = east] (d3) node {$d_3$};
\draw[anchor = south] (E3) node {$e_3$};
\draw[anchor = north] (A4) node {$a_4$};
\draw[anchor = east] (b4) node {$b_4$};
\draw[anchor = east] (c4) node {$c_4$};
\draw[anchor = east] (d4) node {$d_4$};
\draw[anchor = south] (E4) node {$e_4$};
\draw[anchor = north] (A5) node {$a_5$};
\draw[anchor = east] (b5) node {$b_5$};
\draw[anchor = east] (c5) node {$c_5$};
\draw[anchor = east] (d5) node {$d_5$};
\draw[anchor = south] (E5) node {$e_5$};
\draw[anchor = north] (A6) node {$a_6$};
\draw[anchor = east] (b6) node {$b_6$};
\draw[anchor = east] (c6) node {$c_6$};
\draw[anchor = east] (d6) node {$d_6$};
\draw[anchor = south] (E6) node {$e_6$};
\draw[anchor = north] (A7) node {$a_7$};
\draw[anchor = east] (b7) node {$b_7$};
\draw[anchor = east] (c7) node {$c_7$};
\draw[anchor = east] (d7) node {$d_7$};
\draw[anchor = south] (E7) node {$e_7$};
\draw[anchor = north] (A8) node {$a_8$};
\draw[anchor = east] (b8) node {$b_8$};
\draw[anchor = east] (c8) node {$c_8$};
\draw[anchor = east] (d8) node {$d_8$};
\draw[anchor = south] (E8) node {$e_8$};
\end{tikzpicture}
\caption{The set $S$ in the graph $G = P_5 \cp  P_8$ in the proof of Claim~\ref{c:claim9}}
\label{fig:grid9}
\end{center}
\end{figure}
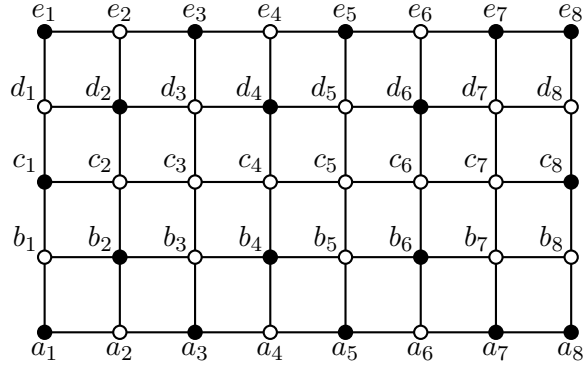

Continuing in this way, the above pattern repeats itself, that is, $S \cap V_i = \{b_i,d_i\}$ for $i$ even and $2 \le i \le m-2$ and $S \cap V_i = \{a_i,e_i\}$ for $i$ odd and $3 \le i \le m-1$. The set $S$ is now fully determined. For example, when $m = 8$ the set $S$ is illustrated in Figure~\ref{fig:grid9}. However, as before the subgraph $G[Z]$ is a path joining two boundary vertices of $G$ and is therefore a $3$-forbidden subgraph of $G$. However, this subgraph contains no vertex of $S$, a contradiction. We deduce, therefore, that our supposition that $m(G,3) = 2m + 2$ is incorrect. Hence, $m(G,3) \ge 2m + 3$. This completes the proof of Claim~\ref{c:claim9}.~\smallqed

\begin{claim}
\label{c:claim10}
If $m$ is even, then $m(G,3) = 2m+3$.
\end{claim}
\proof Suppose that $m$ is even. Let
\[
S_\even = (A_\odd \cup \{a_m\}) \cup (B_\even \setminus \{b_m\}) \cup \{c_1,c_{m-1},c_m\} \cup (D_\even \setminus \{d_m\}) \cup (E_\odd \cup \{a_m\}).
\]

\begin{figure}[htb]
\begin{center}
\begin{tikzpicture}[scale=1,style=thick,x=1cm,y=1cm]
\def\vr{2.5pt} 
\path (0,0) coordinate (A1);
\path (0,1) coordinate (B1);
\path (0.05,1.25) coordinate (b1);
\path (0,2) coordinate (C1);
\path (0.05,2.25) coordinate (c1);
\path (0,3) coordinate (D1);
\path (0.05,3.25) coordinate (d1);
\path (0,4) coordinate (E1);
\path (1,0) coordinate (A2);
\path (1,1) coordinate (B2);
\path (1.05,1.25) coordinate (b2);
\path (1,2) coordinate (C2);
\path (1.05,2.25) coordinate (c2);
\path (1,3) coordinate (D2);
\path (1.05,3.25) coordinate (d2);
\path (1,4) coordinate (E2);
\path (2,0) coordinate (A3);
\path (2,1) coordinate (B3);
\path (2.05,1.25) coordinate (b3);
\path (2,2) coordinate (C3);
\path (2.05,2.25) coordinate (c3);
\path (2,3) coordinate (D3);
\path (2.05,3.25) coordinate (d3);
\path (2,4) coordinate (E3);
\path (3,0) coordinate (A4);
\path (3,1) coordinate (B4);
\path (3.05,1.25) coordinate (b4);
\path (3,2) coordinate (C4);
\path (3.05,2.25) coordinate (c4);
\path (3,3) coordinate (D4);
\path (3.05,3.25) coordinate (d4);
\path (3,4) coordinate (E4);
\path (4,0) coordinate (A5);
\path (4,1) coordinate (B5);
\path (4.05,1.25) coordinate (b5);
\path (4,2) coordinate (C5);
\path (4.05,2.25) coordinate (c5);
\path (4,3) coordinate (D5);
\path (4.05,3.25) coordinate (d5);
\path (4,4) coordinate (E5);
\path (5,0) coordinate (A6);
\path (5,1) coordinate (B6);
\path (5.05,1.25) coordinate (b6);
\path (5,2) coordinate (C6);
\path (5.05,2.25) coordinate (c6);
\path (5,3) coordinate (D6);
\path (5.05,3.25) coordinate (d6);
\path (5,4) coordinate (E6);
\path (6,0) coordinate (A7);
\path (6,1) coordinate (B7);
\path (6.05,1.25) coordinate (b7);
\path (6,2) coordinate (C7);
\path (6.05,2.25) coordinate (c7);
\path (6,3) coordinate (D7);
\path (6.05,3.25) coordinate (d7);
\path (6,4) coordinate (E7);
\path (7,0) coordinate (A8);
\path (7,1) coordinate (B8);
\path (7.05,1.25) coordinate (b8);
\path (7,2) coordinate (C8);
\path (7.05,2.25) coordinate (c8);
\path (7,3) coordinate (D8);
\path (7.05,3.25) coordinate (d8);
\path (7,4) coordinate (E8);
\draw (A1)--(B1)--(C1)--(D1)--(E1);
\draw (A2)--(B2)--(C2)--(D2)--(E2);
\draw (A3)--(B3)--(C3)--(D3)--(E3);
\draw (A4)--(B4)--(C4)--(D4)--(E4);
\draw (A5)--(B5)--(C5)--(D5)--(E5);
\draw (A6)--(B6)--(C6)--(D6)--(E6);
\draw (A7)--(B7)--(C7)--(D7)--(E7);
\draw (A8)--(B8)--(C8)--(D8)--(E8);
\draw (A1)--(A2)--(A3)--(A4)--(A5)--(A6)--(A7)--(A8);
\draw (B1)--(B2)--(B3)--(B4)--(B5)--(B6)--(B7)--(B8);
\draw (C1)--(C2)--(C3)--(C4)--(C5)--(C6)--(C7)--(C8);
\draw (D1)--(D2)--(D3)--(D4)--(D5)--(D6)--(D7)--(D8);
\draw (E1)--(E2)--(E3)--(E4)--(E5)--(E6)--(E7)--(E8);
\draw (A1) [fill=black] circle (\vr);
\draw (B1) [fill=white] circle (\vr);
\draw (C1) [fill=black] circle (\vr);
\draw (D1) [fill=white] circle (\vr);
\draw (E1) [fill=black] circle (\vr);
\draw (A2) [fill=white] circle (\vr);
\draw (B2) [fill=black] circle (\vr);
\draw (C2) [fill=white] circle (\vr);
\draw (D2) [fill=black] circle (\vr);
\draw (E2) [fill=white] circle (\vr);
\draw (A3) [fill=black] circle (\vr);
\draw (B3) [fill=white] circle (\vr);
\draw (C3) [fill=white] circle (\vr);
\draw (D3) [fill=white] circle (\vr);
\draw (E3) [fill=black] circle (\vr);
\draw (A4) [fill=white] circle (\vr);
\draw (B4) [fill=black] circle (\vr);
\draw (C4) [fill=white] circle (\vr);
\draw (D4) [fill=black] circle (\vr);
\draw (E4) [fill=white] circle (\vr);
\draw (A5) [fill=black] circle (\vr);
\draw (B5) [fill=white] circle (\vr);
\draw (C5) [fill=white] circle (\vr);
\draw (D5) [fill=white] circle (\vr);
\draw (E5) [fill=black] circle (\vr);
\draw (A6) [fill=white] circle (\vr);
\draw (B6) [fill=black] circle (\vr);
\draw (C6) [fill=white] circle (\vr);
\draw (D6) [fill=black] circle (\vr);
\draw (E6) [fill=white] circle (\vr);
\draw (A7) [fill=black] circle (\vr);
\draw (B7) [fill=white] circle (\vr);
\draw (C7) [fill=black] circle (\vr);
\draw (D7) [fill=white] circle (\vr);
\draw (E7) [fill=black] circle (\vr);
\draw (A8) [fill=black] circle (\vr);
\draw (B8) [fill=white] circle (\vr);
\draw (C8) [fill=black] circle (\vr);
\draw (D8) [fill=white] circle (\vr);
\draw (E8) [fill=black] circle (\vr);
\draw[anchor = north] (A1) node {$a_1$};
\draw[anchor = east] (b1) node {$b_1$};
\draw[anchor = east] (c1) node {$c_1$};
\draw[anchor = east] (d1) node {$d_1$};
\draw[anchor = south] (E1) node {$e_1$};
\draw[anchor = north] (A2) node {$a_2$};
\draw[anchor = east] (b2) node {$b_2$};
\draw[anchor = east] (c2) node {$c_2$};
\draw[anchor = east] (d2) node {$d_2$};
\draw[anchor = south] (E2) node {$e_2$};
\draw[anchor = north] (A3) node {$a_3$};
\draw[anchor = east] (b3) node {$b_3$};
\draw[anchor = east] (c3) node {$c_3$};
\draw[anchor = east] (d3) node {$d_3$};
\draw[anchor = south] (E3) node {$e_3$};
\draw[anchor = north] (A4) node {$a_4$};
\draw[anchor = east] (b4) node {$b_4$};
\draw[anchor = east] (c4) node {$c_4$};
\draw[anchor = east] (d4) node {$d_4$};
\draw[anchor = south] (E4) node {$e_4$};
\draw[anchor = north] (A5) node {$a_5$};
\draw[anchor = east] (b5) node {$b_5$};
\draw[anchor = east] (c5) node {$c_5$};
\draw[anchor = east] (d5) node {$d_5$};
\draw[anchor = south] (E5) node {$e_5$};
\draw[anchor = north] (A6) node {$a_6$};
\draw[anchor = east] (b6) node {$b_6$};
\draw[anchor = east] (c6) node {$c_6$};
\draw[anchor = east] (d6) node {$d_6$};
\draw[anchor = south] (E6) node {$e_6$};
\draw[anchor = north] (A7) node {$a_7$};
\draw[anchor = east] (b7) node {$b_7$};
\draw[anchor = east] (c7) node {$c_7$};
\draw[anchor = east] (d7) node {$d_7$};
\draw[anchor = south] (E7) node {$e_7$};
\draw[anchor = north] (A8) node {$a_8$};
\draw[anchor = east] (b8) node {$b_8$};
\draw[anchor = east] (c8) node {$c_8$};
\draw[anchor = east] (d8) node {$d_8$};
\draw[anchor = south] (E8) node {$e_8$};
\end{tikzpicture}
\caption{The set $S_\even$ in the graph $G = P_5 \cp P_8$}
\label{fig:grid10}
\end{center}
\end{figure}
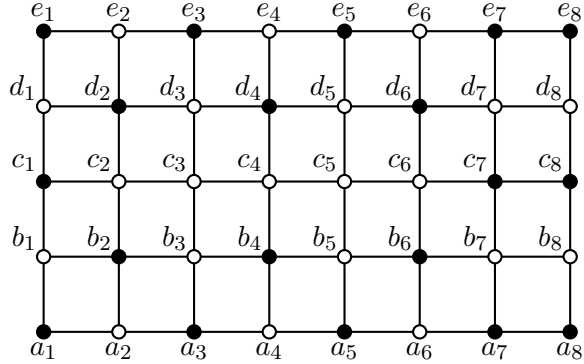

For example, when $m = 8$ the set $S_\even$ is illustrated in Figure~\ref{fig:grid10}.

We show that $S_\even$ is a $3$-percolating set of $G$. The vertices $b_1$ and $d_1$ both have three infected neighbors, and so become infected during the percolation process starting with the initial set $S_\even$. Hence, all vertices in the first column $V_1$ are infected. Every vertex in $A \cup E$ is in the set $S_\even$ or has three infected neighbors, and so become infected during the percolation process. Hence, all vertices in $A \cup E$ are infected. Thereafter, all vertices in $B$ become infected by considering the vertices sequentially (that is, the vertex $b_1$ is first infected, followed by $b_3, b_5, \ldots, b_{m-1}$, and finally $b_m$ is infected. Identical argument show that all vertices in $D$ become infected. Hence all vertices in $B \cup D$ are infected. Thereafter, all vertices in $C$ become infected by considering the vertices sequentially $c_2,c_2,\ldots,c_{m-2}$. Thus, all vertices in $V(G)$ become infected, implying that
\[
\begin{array}{lcl}
m(G,3) \le |S_\even| & = & (|A_\odd|+1) + (|B_\even|-1) + |\{c_1,c_{m-1},c_m\}|  \1 \\
& & \hspace*{0.5cm} + (|D_\even|-1) + (|E_\odd|+1) \1 \\
& = & |A_\odd| + |B_\even| + 3 + |D_\even| + |E_\odd| \1 \\
& = & \frac{m}{2} + \frac{m}{2} + 3 + \frac{m}{2} + \frac{m}{2}  \1 \\
& = & \displaystyle{ 2m+3}.
\end{array}
\]
Hence, $m(G,3) \le 2m+3$. By Claim~\ref{c:claim9}, $m(G,3) \ge 2m+3$. Consequently, $m(G,3) = 2m+3$.~\smallqed

\medskip
The proof of Theorem~\ref{thm:5xm} now follows from Claim~\ref{c:claim8} and~\ref{c:claim10}.~\QED

\subsection{$3$-Bootstrap percolation in $4 \times m$ grids}
\label{S:4xm}

In this section, we show that the $3$-percolation number of a $4 \times m$ grid for all $m \ge 4$ takes on one of two possible values. We first prove a lower bound on the $3$-percolation number of a $4 \times m$ grid.

\begin{thm}
\label{thm:4xm-lower}
For $m \ge 4$, if $G = P_4 \cp P_m$, then $m(G,3) \ge \lfloor \frac{5(m+1)}{3} \rfloor +1$.
\end{thm}
\proof For $m \ge 4$, let $G$ be the grid $P_4 \cp P_m$ with
\[
V(G) = \bigcup_{i=1}^{m} \{a_i,b_i,c_i,d_i\},
\]
where the path $a_ib_ic_id_i$ is a $P_5$-fiber in $G$ for $i \in [m]$, and where $a_1a_2 \ldots a_m$, $b_1b_2 \ldots b_m$, $c_1c_2 \ldots c_m$, and $d_1d_2 \ldots d_m$ are $P_m$-fibers in $G$. For example, when $m = 6$ the grid $G = P_4 \cp P_m$ is illustrated in Figure~\ref{fig:grid6}.

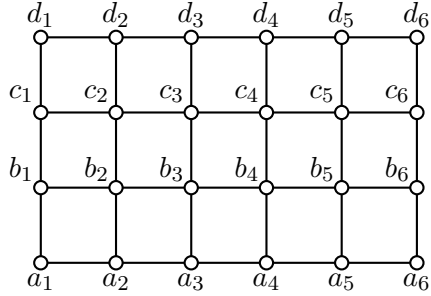
\begin{figure}[htb]
\begin{center}
\begin{tikzpicture}[scale=1,style=thick,x=1cm,y=1cm]
\def\vr{2.5pt} 
\path (0,0) coordinate (A1);
\path (0,1) coordinate (B1);
\path (0.05,1.25) coordinate (b1);
\path (0,2) coordinate (C1);
\path (0.05,2.25) coordinate (c1);
\path (0,3) coordinate (D1);
\path (0.05,3.25) coordinate (d1);
\path (1,0) coordinate (A2);
\path (1,1) coordinate (B2);
\path (1.05,1.25) coordinate (b2);
\path (1,2) coordinate (C2);
\path (1.05,2.25) coordinate (c2);
\path (1,3) coordinate (D2);
\path (1.05,3.25) coordinate (d2);
\path (2,0) coordinate (A3);
\path (2,1) coordinate (B3);
\path (2.05,1.25) coordinate (b3);
\path (2,2) coordinate (C3);
\path (2.05,2.25) coordinate (c3);
\path (2,3) coordinate (D3);
\path (2.05,3.25) coordinate (d3);
\path (3,0) coordinate (A4);
\path (3,1) coordinate (B4);
\path (3.05,1.25) coordinate (b4);
\path (3,2) coordinate (C4);
\path (3.05,2.25) coordinate (c4);
\path (3,3) coordinate (D4);
\path (3.05,3.25) coordinate (d4);
\path (4,0) coordinate (A5);
\path (4,1) coordinate (B5);
\path (4.05,1.25) coordinate (b5);
\path (4,2) coordinate (C5);
\path (4.05,2.25) coordinate (c5);
\path (4,3) coordinate (D5);
\path (4.05,3.25) coordinate (d5);
\path (5,0) coordinate (A6);
\path (5,1) coordinate (B6);
\path (5.05,1.25) coordinate (b6);
\path (5,2) coordinate (C6);
\path (5.05,2.25) coordinate (c6);
\path (5,3) coordinate (D6);
\path (5.05,3.25) coordinate (d6);
\path (6,0) coordinate (A7);
\path (6,1) coordinate (B7);
\path (6.05,1.25) coordinate (b7);
\path (6,2) coordinate (C7);
\path (6.05,2.25) coordinate (c7);
\path (6,3) coordinate (D7);
\path (6.05,3.25) coordinate (d7);
\draw (A1)--(B1)--(C1)--(D1);
\draw (A2)--(B2)--(C2)--(D2);
\draw (A3)--(B3)--(C3)--(D3);
\draw (A4)--(B4)--(C4)--(D4);
\draw (A5)--(B5)--(C5)--(D5);
\draw (A6)--(B6)--(C6)--(D6);
\draw (A1)--(A2)--(A3)--(A4)--(A5)--(A6);
\draw (B1)--(B2)--(B3)--(B4)--(B5)--(B6);
\draw (C1)--(C2)--(C3)--(C4)--(C5)--(C6);
\draw (D1)--(D2)--(D3)--(D4)--(D5)--(D6);
\draw (A1) [fill=white] circle (\vr);
\draw (B1) [fill=white] circle (\vr);
\draw (C1) [fill=white] circle (\vr);
\draw (D1) [fill=white] circle (\vr);
\draw (A2) [fill=white] circle (\vr);
\draw (B2) [fill=white] circle (\vr);
\draw (C2) [fill=white] circle (\vr);
\draw (D2) [fill=white] circle (\vr);
\draw (A3) [fill=white] circle (\vr);
\draw (B3) [fill=white] circle (\vr);
\draw (C3) [fill=white] circle (\vr);
\draw (D3) [fill=white] circle (\vr);
\draw (A4) [fill=white] circle (\vr);
\draw (B4) [fill=white] circle (\vr);
\draw (C4) [fill=white] circle (\vr);
\draw (D4) [fill=white] circle (\vr);
\draw (A5) [fill=white] circle (\vr);
\draw (B5) [fill=white] circle (\vr);
\draw (C5) [fill=white] circle (\vr);
\draw (D5) [fill=white] circle (\vr);
\draw (A6) [fill=white] circle (\vr);
\draw (B6) [fill=white] circle (\vr);
\draw (C6) [fill=white] circle (\vr);
\draw (D6) [fill=white] circle (\vr);
%
\draw[anchor = north] (A1) node {$a_1$};
\draw[anchor = east] (b1) node {$b_1$};
\draw[anchor = east] (c1) node {$c_1$};
\draw[anchor = south] (D1) node {$d_1$};
\draw[anchor = north] (A2) node {$a_2$};
\draw[anchor = east] (b2) node {$b_2$};
\draw[anchor = east] (c2) node {$c_2$};
\draw[anchor = south] (D2) node {$d_2$};
\draw[anchor = north] (A3) node {$a_3$};
\draw[anchor = east] (b3) node {$b_3$};
\draw[anchor = east] (c3) node {$c_3$};
\draw[anchor = south] (D3) node {$d_3$};
\draw[anchor = north] (A4) node {$a_4$};
\draw[anchor = east] (b4) node {$b_4$};
\draw[anchor = east] (c4) node {$c_4$};
\draw[anchor = south] (D4) node {$d_4$};
\draw[anchor = north] (A5) node {$a_5$};
\draw[anchor = east] (b5) node {$b_5$};
\draw[anchor = east] (c5) node {$c_5$};
\draw[anchor = south] (D5) node {$d_5$};
\draw[anchor = north] (A6) node {$a_6$};
\draw[anchor = east] (b6) node {$b_6$};
\draw[anchor = east] (c6) node {$c_6$};
\draw[anchor = south] (D6) node {$d_6$};
\end{tikzpicture}
\caption{The graph $G = P_4 \cp  P_6$}
\label{fig:grid11}
\end{center}
\end{figure}

For $i \in [m]$, let $V_i = \{a_i,b_i,c_i,d_i\}$ and let
\[
V_{\le i }  = \bigcup_{j=1}^i V_i \hspace*{0.5cm} \mbox{and} \hspace*{0.5cm} V_{\ge i }  = \bigcup_{j=i}^m V_i.
\]

Thus, $V(G) = V_{\le m} = V_{\ge 1}$. Let $A = \{a_1,a_2,\ldots,a_m\}$, $B = \{b_1,b_2,\ldots,b_m\}$, $C = \{c_1,c_2,\ldots,c_m\}$, and $D = \{d_1,d_2,\ldots,d_m\}$. In what follows, let $S$ be a minimum $3$-percolating set of $G$ that does not contain three consecutive boundary vertices of $G$. We note that such a set $S$ exists by Lemma~\ref{lem:three boundary vertices}.

\begin{claim}
\label{c:claim11}
The following properties hold. \\ [-26pt]
\begin{enumerate}
\item[{\rm (a)}] $|S \cap V_1| \ge 3$ and $|S \cap V_m| \ge 3$.
\item[{\rm (b)}] $|S \cap V_i| \ge 2$ for all $i$ where $2 \le i \le m-1$;
\item[{\rm (c)}] $|S \cap (V_i \cup V_{i+1})| \ge 3$ for all $i$ where $2 \le i \le m-1$.
\item[{\rm (d)}] $|S \cap (V_i \cup V_{i+1} \cup V_{i+2})| \ge 5$ for all $i$ where $2 \le i \le m-2$.
\end{enumerate}
\end{claim}
\proof Since the vertices $a_1$ and $d_1$ both have degree~$2$ in $G$, we note that $\{a_1,d_1\} \subset S$. Since the set $S$ contains at least one of every two adjacent boundary vertices, we note that $|S \cap \{b_1,c_1\}| \ge 1$, implying that $|S \cap V_1| \ge 3$. By symmetry, $|S \cap V_m| \ge 3$. Thus, Property~(a) holds. Property~(b) follows from Corollary~\ref{cor:subgraph condition_1}(b).

To prove Property~(c), consider the set $S \cap (V_i \cup V_{i+1})$ for some $i$ where $2 \le i \le m-1$. By Corollary~\ref{cor:subgraph condition_1}(a), $|S \cap \{a_i,a_{i+1}\}| \ge 1$ and $|S \cap \{d_i,d_{i+1}\}| \ge 1$. If $Q_1 = G[\lbrace b_i,b_{i+1},c_i,c_{i+1} \rbrace]$, then $Q_1 = C_4$, and so the subgraphs $Q_1$ is a $3$-forbidden subgraph of $G$, and so $S$ must contain at least one vertex from each of $Q_1$. These observations imply that $|S \cap (V_i \cup V_{i+1})| \ge 3$, and so Property~(c) holds.

To prove Property~(d), consider the set $S \cap (V_i \cup V_{i+1} \cup V_{i+2})$ for some $i$ where $2 \le i \le m-2$. For notation convenience, we may assume that $i = 2$, that is, we consider the set $S \cap (V_2 \cup V_3 \cup V_4)$. If $S$ contains at least three boundary vertices in $V_2 \cup V_3 \cup V_4$, then the desired result is immediate. Hence, we may assume that $|S \cap \{a_2,a_3,a_4,d_2,d_3,d_4\}| \le 4$, for otherwise the desired lower bound holds. Since $S$ contains no three consecutive boundaries, we note that  $|S \cap \{a_2,a_3,a_4\}| \le 2$ and $|S \cap \{d_2,d_3,d_4\}| \le 2$.

Suppose that $S$ contains four boundary vertices in $V_2 \cup V_3 \cup V_4$, implying that $\{a_2,a_4,d_2,d_4\} \subset S$. By Property~(b) we have $|S \cap V_3| \ge 1$, we infer in this case that $|S \cap (V_2 \cup V_{2} \cup V_{4})| \ge 5$. Hence, we may assume that $S$ contains at most three boundary vertices in $V_2 \cup V_3 \cup V_4$, for otherwise the desired lower bound holds.

Suppose that $S$ contains exactly three boundary vertices in $V_2 \cup V_3 \cup V_4$. By symmetry, we may assume that $\{a_2,a_4,d_2\} \subset S$ or $\{a_2,a_4,d_3\} \subset S$. Suppose that $\{a_2,a_4,d_2\} \subset S$. Let $Q_2 = G[\lbrace b_2,b_3,c_2,c_3 \rbrace]$, $Q_3 = G[\lbrace c_3,c_4,d_3,d_4 \rbrace]$, and $Q_4 = G[\lbrace a_3,b_3,b_4,c_4,d_4 \rbrace]$. Since $Q_2 = Q_3 = C_4$ and since $Q_3$ is a path joining two boundary vertices of $G$, the subgraphs $Q_2$, $Q_3$ and $Q_4$ are all $3$-forbidden subgraphs of $G$, and so $S$ must contain at least one vertex from each of $Q_2$, $Q_3$ and $Q_4$. At least two vertices in $S$ are needed for this purpose, implying that $|S \cap (V_2 \cup V_3 \cup V_4)| \ge 5$, as desired. Suppose next that $\{a_2,a_4,d_3\} \subset S$. Let $Q_5 = G[\lbrace d_2,c_2,c_3,c_4,d_4 \rbrace]$, $Q_6 = G[\lbrace d_2,c_2,b_2,b_3,a_3 \rbrace]$, and $Q_7 = G[\lbrace b_3,b_4,c_3,c_4 \rbrace]$. Since $Q_5$ and $Q_6$ are path joining two boundary vertices of $G$ and since $Q_7 = C_4$, the subgraphs $Q_5$, $Q_6$ and $Q_7$ are all $3$-forbidden subgraphs of $G$, and so $S$ must contain at least one vertex from each of $Q_5$, $Q_6$ and $Q_7$. At least two vertices in $S$ are needed for this purpose, implying once again that $|S \cap (V_2 \cup V_3 \cup V_4)| \ge 5$, as desired.

Hence, we may assume that $S$ contains at most two boundary vertices in $V_2 \cup V_3 \cup V_4$, for otherwise the desired lower bound holds. Since $S$ contains at least one vertex among every two adjacent boundary vertices, this implies that $a_3$ and $d_3$ are the two boundary vertices in $S$. We note that $|S \cap V_2| \ge 1$ and $|S \cap V_4| \ge 1$. Suppose that $|S \cap V_2| = 1$ and $|S \cap V_4| = 1$, implying by our earlier assumptions that $|S \cap \{b_2,c_2\}| = 1$ and $|S \cap \{b_4,c_4\}| = 1$. By symmetry, we may assume that $\{b_2,c_4\} \subset S$ or $\{b_2,b_4\} \subset S$. If $\{b_2,c_4\} \subset S$, then $G[\lbrace d_2,c_2,c_3,b_3,b_4,a_4 \rbrace]$ is a path joining two boundary vertices that contains no vertex of $S$, and if $\{b_2,b_4\} \subset S$, then $G[\lbrace d_2,c_2,c_3,c_4,d_4 \rbrace]$ is a path joining two boundary vertices  that contains no vertex of $S$. Both cases produce a contradiction. We deduce, therefore, that $|S \cap V_2| \ge 2$ or $|S \cap V_4| \ge 2$, implying that $|S \cap (V_2 \cup V_3 \cup V_4)| \ge 5$, as desired. This completes the proof of Property~(d).~\smallqed

\medskip
We now return to the proof of Theorem~\ref{thm:4xm-lower} and calculate the lower bound on $m(G,3)$.

\begin{claim}
\label{c:claim12}
If $m \equiv 0 \, (\modo \, 3)$, then $m(G,3) \ge \lfloor \frac{5(m+1)}{3} \rfloor + 1$.
\end{claim}
\proof Suppose that $m \equiv 0 \, (\modo \, 3)$. By Claim~\ref{c:claim11} we have
\[
\begin{array}{lcl}
m(G,3) = |S| & = & \displaystyle{ |S \cap V_1| + |S \cap V_2| + |S \cap V_m| + \sum_{i=3}^{m-1} |S \cap V_i| }  \1 \\
& \ge & \displaystyle{ 3 + 1 + 3 + \sum_{i=1}^{\frac{m-3}{3}} |S \cap (V_{3i} \cup V_{3i+1} \cup V_{3i+2} )|     }   \1 \\
& \ge & 3 + 1 + 3 + \frac{m-3}{3} \times 5   \1 \\
& = & \frac{5}{3}(m+1) + \frac{1}{3}  \1 \\
& = & \lfloor \frac{5(m+1)}{3} \rfloor + 1,
\end{array}
\]
noting that in this case $m \equiv 0 \, (\modo \, 3)$.~\smallqed

\begin{claim}
\label{c:claim13}
If $m \equiv 1 \, (\modo \, 3)$, then $m(G,3) \ge \lfloor \frac{5(m+1)}{3} \rfloor + 1$.
\end{claim}
\proof Suppose that $m \equiv 1 \, (\modo \, 3)$. By Claim~\ref{c:claim11} we have
\[
\begin{array}{lcl}
m(G,3) = |S| & = & \displaystyle{ |S \cap V_1| + |S \cap (V_2 \cup V_3)| + |S \cap V_m| + \sum_{i=4}^{m-1} |S \cap V_i| }  \1 \\
& \ge & \displaystyle{ 3 + 3 + 3 + \sum_{i=1}^{\frac{m-4}{3}} |S \cap (V_{3i+1} \cup V_{3i+2} \cup V_{3i+3} )|     }   \1 \\
& \ge & 3 + 3 + 3 + \frac{m-4}{3} \times 5   \1 \\
& = & \frac{5}{3}(m+1) + \frac{2}{3} \2 \\
& = & \lfloor \frac{5(m+1)}{3} \rfloor + 1, \1
\end{array}
\]
noting that in this case $m \equiv 1 \, (\modo \, 3)$.~\smallqed

\newpage
\begin{claim}
\label{c:claim14}
If $m \equiv 2 \, (\modo \, 3)$, then $m(G,3) \ge \lfloor \frac{5(m+1)}{3} \rfloor + 1$.
\end{claim}
\proof Suppose that $m \equiv 2 \, (\modo \, 3)$. By Claim~\ref{c:claim11} we have
\[
\begin{array}{lcl}
m(G,3) = |S| & = & \displaystyle{ |S \cap V_1| + |S \cap V_m| + \sum_{i=2}^{m-1} |S \cap V_i| }  \1 \\
& \ge & \displaystyle{ 3 +3 + \sum_{i=1}^{\frac{m-2}{3}} |S \cap (V_{3i-1} \cup V_{3i} \cup V_{3i+1} )|     }   \1 \\
& \ge & 3 +3 + \frac{m-2}{3} \times 5   \1 \\
& = & \frac{5}{3}(m+1) + 1  \2 \\
& = & \lfloor \frac{5(m+1)}{3} \rfloor + 1, \1
\end{array}
\]
noting that in this case $m \equiv 2 \, (\modo \, 3)$.~\smallqed

\medskip
By Claims~\ref{c:claim12},~\ref{c:claim13}, and~\ref{c:claim14}, we have $m(G,3) \ge \lfloor \frac{5(m+1)}{3} \rfloor + 1$. This completes the proof of Theorem~\ref{thm:4xm-lower}.~\QED

We establish next upper bounds on the $3$-percolation number of $4 \times m$ grids for all $m \ge 4$.

\begin{thm}
\label{thm:4xm-upper}
For $m \ge 4$, if $G = P_4 \cp P_m$, then
\[
m(G,3) \le
\begin{cases}
        \lfloor \frac{5(m+1)}{3} \rfloor + 1; &   m \in \{5,7,11\} \2 \\
        \lfloor \frac{5(m+1)}{3} \rfloor + 2; &   \text{ otherwise}.
    \end{cases}
\]
\end{thm}
\proof For $m \ge 4$, let $G_m$ be the grid $P_4 \cp P_m$ where we follow the notation in the proof of Theorem~\ref{thm:4xm-lower}. The sets shown in Figure~\ref{fig:S for m = 5,7,11}(a),~\ref{fig:S for m = 5,7,11}(b), and~\ref{fig:S for m = 5,7,11}(c) are $3$-percolating sets of $G_5$, $G_7$, and $G_{11}$, respectively, of cardinalities~$11$,~$14$ and~$21$, respectively, implying that
\[
m(G_m,3) \le \left\lfloor \frac{5(m+1)}{3} \right\rfloor + 1
\]
for $m \in \{5,7,11\}$. Hence in what follows, we may assume that $m \notin \{5,7,11\}$, for otherwise the desired upper bound follows. For $i \in \{2,\ldots,m-3\}$, let
\[
X_i = \{a_i,c_i,b_{i+1},d_{i+1},a_{i+2}\} \hspace*{0.5cm} \mbox{and} \hspace*{0.5cm}
Y_i = \{b_i,d_i,a_{i+1},c_{i+1},d_{i+2}\}.
\]

\begin{figure}[htb]
\begin{center}
\begin{tikzpicture}[scale=0.9,style=thick,x=1cm,y=1cm]
\def\vr{2.5pt} 


\path (0,0) coordinate (a1);
\path (0,1) coordinate (b1);
\path (0,2) coordinate (c1);
\path (0,3) coordinate (d1);

\path (1,0) coordinate (a2);
\path (1,1) coordinate (b2);
\path (1,2) coordinate (c2);
\path (1,3) coordinate (d2);

\path (2,0) coordinate (a3);
\path (2,1) coordinate (b3);
\path (2,2) coordinate (c3);
\path (2,3) coordinate (d3);

\path (3,0) coordinate (a4);
\path (3,1) coordinate (b4);
\path (3,2) coordinate (c4);
\path (3,3) coordinate (d4);

\path (4,0) coordinate (a5);
\path (4,1) coordinate (b5);
\path (4,2) coordinate (c5);
\path (4,3) coordinate (d5);

\path (5,0) coordinate (a6);
\path (5,1) coordinate (b6);
\path (5,2) coordinate (c6);
\path (5,3) coordinate (d6);

\path (6,0) coordinate (a7);
\path (6,1) coordinate (b7);
\path (6,2) coordinate (c7);
\path (6,3) coordinate (d7);

\path (7,0) coordinate (a8);
\path (7,1) coordinate (b8);
\path (7,2) coordinate (c8);
\path (7,3) coordinate (d8);

\path (8,0) coordinate (a9);
\path (8,1) coordinate (b9);
\path (8,2) coordinate (c9);
\path (8,3) coordinate (d9);

\path (9,0) coordinate (a10);
\path (9,1) coordinate (b10);
\path (9,2) coordinate (c10);
\path (9,3) coordinate (d10);

\path (10,0) coordinate (a11);
\path (10,1) coordinate (b11);
\path (10,2) coordinate (c11);
\path (10,3) coordinate (d11);

\draw (a1)--(b1)--(c1)--(d1);
\draw (a2)--(b2)--(c2)--(d2);
\draw (a3)--(b3)--(c3)--(d3);
\draw (a4)--(b4)--(c4)--(d4);
\draw (a5)--(b5)--(c5)--(d5);
\draw (a6)--(b6)--(c6)--(d6);
\draw (a7)--(b7)--(c7)--(d7);
\draw (a8)--(b8)--(c8)--(d8);
\draw (a9)--(b9)--(c9)--(d9);
\draw (a10)--(b10)--(c10)--(d10);
\draw (a11)--(b11)--(c11)--(d11);

\draw (a1)--(a2)--(a3)--(a4)--(a5)--(a6)--(a7)--(a8)--(a9)--(a10)--(a11);
\draw (b1)--(b2)--(b3)--(b4)--(b5)--(b6)--(b7)--(b8)--(b9)--(b10)--(b11);
\draw (c1)--(c2)--(c3)--(c4)--(c5)--(c6)--(c7)--(c8)--(c9)--(c10)--(c11);
\draw (d1)--(d2)--(d3)--(d4)--(d5)--(d6)--(d7)--(d8)--(d9)--(d10)--(d11);

\draw (a1) [fill=black] circle (\vr);
\draw (b1) [fill=black] circle (\vr);
\draw (c1) [fill=white] circle (\vr);
\draw (d1) [fill=black] circle (\vr);

\draw (a2) [fill=white] circle (\vr);
\draw (b2) [fill=white] circle (\vr);
\draw (c2) [fill=black] circle (\vr);
\draw (d2) [fill=white] circle (\vr);

\draw (a3) [fill=black] circle (\vr);
\draw (b3) [fill=white] circle (\vr);
\draw (c3) [fill=white] circle (\vr);
\draw (d3) [fill=black] circle (\vr);

\draw (a4) [fill=white] circle (\vr);
\draw (b4) [fill=black] circle (\vr);
\draw (c4) [fill=black] circle (\vr);
\draw (d4) [fill=white] circle (\vr);

\draw (a5) [fill=black] circle (\vr);
\draw (b5) [fill=white] circle (\vr);
\draw (c5) [fill=white] circle (\vr);
\draw (d5) [fill=black] circle (\vr);

\draw (a6) [fill=white] circle (\vr);
\draw (b6) [fill=white] circle (\vr);
\draw (c6) [fill=black] circle (\vr);
\draw (d6) [fill=white] circle (\vr);

\draw (a7) [fill=black] circle (\vr);
\draw (b7) [fill=white] circle (\vr);
\draw (c7) [fill=white] circle (\vr);
\draw (d7) [fill=black] circle (\vr);

\draw (a8) [fill=white] circle (\vr);
\draw (b8) [fill=black] circle (\vr);
\draw (c8) [fill=black] circle (\vr);
\draw (d8) [fill=white] circle (\vr);

\draw (a9) [fill=black] circle (\vr);
\draw (b9) [fill=white] circle (\vr);
\draw (c9) [fill=white] circle (\vr);
\draw (d9) [fill=black] circle (\vr);

\draw (a10) [fill=white] circle (\vr);
\draw (b10) [fill=white] circle (\vr);
\draw (c10) [fill=black] circle (\vr);
\draw (d10) [fill=white] circle (\vr);

\draw (a11) [fill=black] circle (\vr);
\draw (b11) [fill=black] circle (\vr);
\draw (c11) [fill=white] circle (\vr);
\draw (d11) [fill=black] circle (\vr);


\path (-1,5) coordinate (a11);
\path (-1,6) coordinate (b11);
\path (-1,7) coordinate (c11);
\path (-1,8) coordinate (d11);

\path (0,5) coordinate (a21);
\path (0,6) coordinate (b21);
\path (0,7) coordinate (c21);
\path (0,8) coordinate (d21);

\path (1,5) coordinate (a31);
\path (1,6) coordinate (b31);
\path (1,7) coordinate (c31);
\path (1,8) coordinate (d31);

\path (2,5) coordinate (a41);
\path (2,6) coordinate (b41);
\path (2,7) coordinate (c41);
\path (2,8) coordinate (d41);

\path (3,5) coordinate (a51);
\path (3,6) coordinate (b51);
\path (3,7) coordinate (c51);
\path (3,8) coordinate (d51);

\draw (a11)--(b11)--(c11)--(d11);
\draw (a21)--(b21)--(c21)--(d21);
\draw (a31)--(b31)--(c31)--(d31);
\draw (a41)--(b41)--(c41)--(d41);
\draw (a51)--(b51)--(c51)--(d51);

\draw (a11)--(a21)--(a31)--(a41)--(a51);
\draw (b11)--(b21)--(b31)--(b41)--(b51);
\draw (c11)--(c21)--(c31)--(c41)--(c51);
\draw (d11)--(d21)--(d31)--(d41)--(d51);

\draw (a11) [fill=black] circle (\vr);
\draw (b11) [fill=black] circle (\vr);
\draw (c11) [fill=white] circle (\vr);
\draw (d11) [fill=black] circle (\vr);

\draw (a21) [fill=white] circle (\vr);
\draw (b21) [fill=white] circle (\vr);
\draw (c21) [fill=black] circle (\vr);
\draw (d21) [fill=white] circle (\vr);

\draw (a31) [fill=black] circle (\vr);
\draw (b31) [fill=white] circle (\vr);
\draw (c31) [fill=white] circle (\vr);
\draw (d31) [fill=black] circle (\vr);

\draw (a41) [fill=white] circle (\vr);
\draw (b41) [fill=black] circle (\vr);
\draw (c41) [fill=black] circle (\vr);
\draw (d41) [fill=white] circle (\vr);

\draw (a51) [fill=black] circle (\vr);
\draw (b51) [fill=black] circle (\vr);
\draw (c51) [fill=white] circle (\vr);
\draw (d51) [fill=black] circle (\vr);


\path (5,5) coordinate (a12);
\path (5,6) coordinate (b12);
\path (5,7) coordinate (c12);
\path (5,8) coordinate (d12);

\path (6,5) coordinate (a22);
\path (6,6) coordinate (b22);
\path (6,7) coordinate (c22);
\path (6,8) coordinate (d22);

\path (7,5) coordinate (a32);
\path (7,6) coordinate (b32);
\path (7,7) coordinate (c32);
\path (7,8) coordinate (d32);

\path (8,5) coordinate (a42);
\path (8,6) coordinate (b42);
\path (8,7) coordinate (c42);
\path (8,8) coordinate (d42);

\path (9,5) coordinate (a52);
\path (9,6) coordinate (b52);
\path (9,7) coordinate (c52);
\path (9,8) coordinate (d52);

\path (10,5) coordinate (a62);
\path (10,6) coordinate (b62);
\path (10,7) coordinate (c62);
\path (10,8) coordinate (d62);

\path (11,5) coordinate (a72);
\path (11,6) coordinate (b72);
\path (11,7) coordinate (c72);
\path (11,8) coordinate (d72);

\draw (a12)--(b12)--(c12)--(d12);
\draw (a22)--(b22)--(c22)--(d22);
\draw (a32)--(b32)--(c32)--(d32);
\draw (a42)--(b42)--(c42)--(d42);
\draw (a52)--(b52)--(c52)--(d52);
\draw (a62)--(b62)--(c62)--(d62);
\draw (a72)--(b72)--(c72)--(d72);

\draw (a12)--(a22)--(a32)--(a42)--(a52)--(a62)--(a72);
\draw (b12)--(b22)--(b32)--(b42)--(b52)--(b62)--(b72);
\draw (c12)--(c22)--(c32)--(c42)--(c52)--(c62)--(c72);
\draw (d12)--(d22)--(d32)--(d42)--(d52)--(d62)--(d72);

\draw (a12) [fill=black] circle (\vr);
\draw (b12) [fill=black] circle (\vr);
\draw (c12) [fill=white] circle (\vr);
\draw (d12) [fill=black] circle (\vr);

\draw (a22) [fill=white] circle (\vr);
\draw (b22) [fill=white] circle (\vr);
\draw (c22) [fill=black] circle (\vr);
\draw (d22) [fill=white] circle (\vr);

\draw (a32) [fill=black] circle (\vr);
\draw (b32) [fill=white] circle (\vr);
\draw (c32) [fill=white] circle (\vr);
\draw (d32) [fill=black] circle (\vr);

\draw (a42) [fill=white] circle (\vr);
\draw (b42) [fill=black] circle (\vr);
\draw (c42) [fill=black] circle (\vr);
\draw (d42) [fill=white] circle (\vr);

\draw (a52) [fill=black] circle (\vr);
\draw (b52) [fill=white] circle (\vr);
\draw (c52) [fill=white] circle (\vr);
\draw (d52) [fill=black] circle (\vr);

\draw (a62) [fill=white] circle (\vr);
\draw (b62) [fill=white] circle (\vr);
\draw (c62) [fill=black] circle (\vr);
\draw (d62) [fill=white] circle (\vr);

\draw (a72) [fill=black] circle (\vr);
\draw (b72) [fill=black] circle (\vr);
\draw (c72) [fill=white] circle (\vr);
\draw (d72) [fill=black] circle (\vr);

\draw (1,4.5) node {{\small (a)}};
\draw (8,4.5) node {{\small (b)}};
\draw (5,-0.5) node {{\small (c)}};
\end{tikzpicture}
\caption{$3$-Percolating sets for $P_4 \cp P_5$, $P_4 \cp P_7$, and $P_4 \cp P_{11}$}
\label{fig:S for m = 5,7,11}
\end{center}
\end{figure}
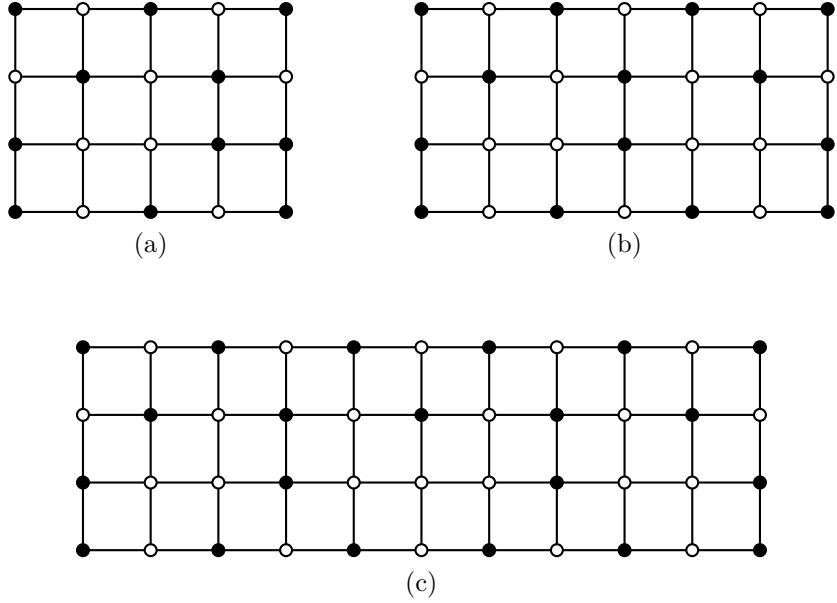

\begin{figure}[htb]
\begin{center}
\begin{tikzpicture}[scale=1,style=thick,x=1cm,y=1cm]
\def\vr{2.5pt} 
\path (0,0) coordinate (A1);
\path (0,1) coordinate (B1);
\path (0.05,1.25) coordinate (b1);
\path (0,2) coordinate (C1);
\path (0.05,2.25) coordinate (c1);
\path (0,3) coordinate (D1);
\path (0.05,3.25) coordinate (d1);
\path (1,0) coordinate (A2);
\path (1,1) coordinate (B2);
\path (1.05,1.25) coordinate (b2);
\path (1,2) coordinate (C2);
\path (1.05,2.25) coordinate (c2);
\path (1,3) coordinate (D2);
\path (1.05,3.25) coordinate (d2);
\path (2,0) coordinate (A3);
\path (2,1) coordinate (B3);
\path (2.05,1.25) coordinate (b3);
\path (2,2) coordinate (C3);
\path (2.05,2.25) coordinate (c3);
\path (2,3) coordinate (D3);
\path (2.05,3.25) coordinate (d3);
%
\draw (A1)--(B1)--(C1)--(D1);
\draw (A2)--(B2)--(C2)--(D2);
\draw (A3)--(B3)--(C3)--(D3);
\draw (A1)--(A2)--(A3);
\draw (B1)--(B2)--(B3);
\draw (C1)--(C2)--(C3);
\draw (D1)--(D2)--(D3);
\draw (A1) [fill=black] circle (\vr);
\draw (B1) [fill=white] circle (\vr);
\draw (C1) [fill=black] circle (\vr);
\draw (D1) [fill=white] circle (\vr);
\draw (A2) [fill=white] circle (\vr);
\draw (B2) [fill=black] circle (\vr);
\draw (C2) [fill=white] circle (\vr);
\draw (D2) [fill=black] circle (\vr);
\draw (A3) [fill=black] circle (\vr);
\draw (B3) [fill=white] circle (\vr);
\draw (C3) [fill=white] circle (\vr);
\draw (D3) [fill=white] circle (\vr);
%
\draw[anchor = north] (A1) node {$a_2$};
\draw[anchor = east] (b1) node {$b_2$};
\draw[anchor = east] (c1) node {$c_2$};
\draw[anchor = south] (D1) node {$d_2$};
\draw[anchor = north] (A2) node {$a_3$};
\draw[anchor = east] (b2) node {$b_3$};
\draw[anchor = east] (c2) node {$c_3$};
\draw[anchor = south] (D2) node {$d_3$};
\draw[anchor = north] (A3) node {$a_4$};
\draw[anchor = east] (b3) node {$b_4$};
\draw[anchor = east] (c3) node {$c_4$};
\draw[anchor = south] (D3) node {$d_4$};
\draw (1,-1) node {{\small (a) $X_2$}};
\path (4,0) coordinate (A1);
\path (4,1) coordinate (B1);
\path (4.05,1.25) coordinate (b1);
\path (4,2) coordinate (C1);
\path (4.05,2.25) coordinate (c1);
\path (4,3) coordinate (D1);
\path (4.05,3.25) coordinate (d1);
\path (5,0) coordinate (A2);
\path (5,1) coordinate (B2);
\path (5.05,1.25) coordinate (b2);
\path (5,2) coordinate (C2);
\path (5.05,2.25) coordinate (c2);
\path (5,3) coordinate (D2);
\path (5.05,3.25) coordinate (d2);
\path (6,0) coordinate (A3);
\path (6,1) coordinate (B3);
\path (6.05,1.25) coordinate (b3);
\path (6,2) coordinate (C3);
\path (6.05,2.25) coordinate (c3);
\path (6,3) coordinate (D3);
\path (6.05,3.25) coordinate (d3);
%
\draw (A1)--(B1)--(C1)--(D1);
\draw (A2)--(B2)--(C2)--(D2);
\draw (A3)--(B3)--(C3)--(D3);
\draw (A1)--(A2)--(A3);
\draw (B1)--(B2)--(B3);
\draw (C1)--(C2)--(C3);
\draw (D1)--(D2)--(D3);
\draw (A1) [fill=white] circle (\vr);
\draw (B1) [fill=black] circle (\vr);
\draw (C1) [fill=white] circle (\vr);
\draw (D1) [fill=black] circle (\vr);
\draw (A2) [fill=black] circle (\vr);
\draw (B2) [fill=white] circle (\vr);
\draw (C2) [fill=black] circle (\vr);
\draw (D2) [fill=white] circle (\vr);
\draw (A3) [fill=white] circle (\vr);
\draw (B3) [fill=white] circle (\vr);
\draw (C3) [fill=white] circle (\vr);
\draw (D3) [fill=black] circle (\vr);
%
\draw[anchor = north] (A1) node {$a_2$};
\draw[anchor = east] (b1) node {$b_2$};
\draw[anchor = east] (c1) node {$c_2$};
\draw[anchor = south] (D1) node {$d_2$};
\draw[anchor = north] (A2) node {$a_3$};
\draw[anchor = east] (b2) node {$b_3$};
\draw[anchor = east] (c2) node {$c_3$};
\draw[anchor = south] (D2) node {$d_3$};
\draw[anchor = north] (A3) node {$a_4$};
\draw[anchor = east] (b3) node {$b_4$};
\draw[anchor = east] (c3) node {$c_4$};
\draw[anchor = south] (D3) node {$d_4$};
\draw (5,-1) node {{\small (b) $Y_2$}};
\path (8,0) coordinate (A1);
\path (8,1) coordinate (B1);
\path (8.05,1.25) coordinate (b1);
\path (8,2) coordinate (C1);
\path (8.05,2.25) coordinate (c1);
\path (8,3) coordinate (D1);
\path (8.05,3.25) coordinate (d1);
\path (9,0) coordinate (A2);
\path (9,1) coordinate (B2);
\path (9.05,1.25) coordinate (b2);
\path (9,2) coordinate (C2);
\path (9.05,2.25) coordinate (c2);
\path (9,3) coordinate (D2);
\path (9.05,3.25) coordinate (d2);
\path (10,0) coordinate (A3);
\path (10,1) coordinate (B3);
\path (10.05,1.25) coordinate (b3);
\path (10,2) coordinate (C3);
\path (10.05,2.25) coordinate (c3);
\path (10,3) coordinate (D3);
\path (10.05,3.25) coordinate (d3);
\path (11,0) coordinate (A4);
\path (11,1) coordinate (B4);
\path (11.05,1.25) coordinate (b4);
\path (11,2) coordinate (C4);
\path (11.05,2.25) coordinate (c4);
\path (11,3) coordinate (D4);
\path (11.05,3.25) coordinate (d4);
\path (12,0) coordinate (A5);
\path (12,1) coordinate (B5);
\path (12.05,1.25) coordinate (b5);
\path (12,2) coordinate (C5);
\path (12.05,2.25) coordinate (c5);
\path (12,3) coordinate (D5);
\path (12.05,3.25) coordinate (d5);
\path (13,0) coordinate (A6);
\path (13,1) coordinate (B6);
\path (13.05,1.25) coordinate (b6);
\path (13,2) coordinate (C6);
\path (13.05,2.25) coordinate (c6);
\path (13,3) coordinate (D6);
\path (13.05,3.25) coordinate (d6);
\path (14,0) coordinate (A7);
\path (14,1) coordinate (B7);
\path (14.05,1.25) coordinate (b7);
\path (14,2) coordinate (C7);
\path (14.05,2.25) coordinate (c7);
\path (14,3) coordinate (D7);
\path (14.05,3.25) coordinate (d7);
\draw (A1)--(B1)--(C1)--(D1);
\draw (A2)--(B2)--(C2)--(D2);
\draw (A3)--(B3)--(C3)--(D3);
\draw (A4)--(B4)--(C4)--(D4);
\draw (A5)--(B5)--(C5)--(D5);
\draw (A6)--(B6)--(C6)--(D6);
\draw (A1)--(A2)--(A3)--(A4)--(A5)--(A6);
\draw (B1)--(B2)--(B3)--(B4)--(B5)--(B6);
\draw (C1)--(C2)--(C3)--(C4)--(C5)--(C6);
\draw (D1)--(D2)--(D3)--(D4)--(D5)--(D6);
\draw (A1) [fill=black] circle (\vr);
\draw (B1) [fill=white] circle (\vr);
\draw (C1) [fill=black] circle (\vr);
\draw (D1) [fill=white] circle (\vr);
\draw (A2) [fill=white] circle (\vr);
\draw (B2) [fill=black] circle (\vr);
\draw (C2) [fill=white] circle (\vr);
\draw (D2) [fill=black] circle (\vr);
\draw (A3) [fill=black] circle (\vr);
\draw (B3) [fill=white] circle (\vr);
\draw (C3) [fill=white] circle (\vr);
\draw (D3) [fill=white] circle (\vr);
\draw (A4) [fill=white] circle (\vr);
\draw (B4) [fill=black] circle (\vr);
\draw (C4) [fill=white] circle (\vr);
\draw (D4) [fill=black] circle (\vr);
\draw (A5) [fill=black] circle (\vr);
\draw (B5) [fill=white] circle (\vr);
\draw (C5) [fill=black] circle (\vr);
\draw (D5) [fill=white] circle (\vr);
\draw (A6) [fill=white] circle (\vr);
\draw (B6) [fill=white] circle (\vr);
\draw (C6) [fill=white] circle (\vr);
\draw (D6) [fill=black] circle (\vr);
%
\draw[anchor = north] (A1) node {$a_2$};
\draw[anchor = east] (b1) node {$b_2$};
\draw[anchor = east] (c1) node {$c_2$};
\draw[anchor = south] (D1) node {$d_2$};
\draw[anchor = north] (A2) node {$a_3$};
\draw[anchor = east] (b2) node {$b_3$};
\draw[anchor = east] (c2) node {$c_3$};
\draw[anchor = south] (D2) node {$d_3$};
\draw[anchor = north] (A3) node {$a_4$};
\draw[anchor = east] (b3) node {$b_4$};
\draw[anchor = east] (c3) node {$c_4$};
\draw[anchor = south] (D3) node {$d_4$};
\draw[anchor = north] (A4) node {$a_5$};
\draw[anchor = east] (b4) node {$b_5$};
\draw[anchor = east] (c4) node {$c_5$};
\draw[anchor = south] (D4) node {$d_5$};
\draw[anchor = north] (A5) node {$a_6$};
\draw[anchor = east] (b5) node {$b_6$};
\draw[anchor = east] (c5) node {$c_6$};
\draw[anchor = south] (D5) node {$d_6$};
\draw[anchor = north] (A6) node {$a_7$};
\draw[anchor = east] (b6) node {$b_7$};
\draw[anchor = east] (c6) node {$c_7$};
\draw[anchor = south] (D6) node {$d_7$};
\draw (11,-1) node {{\small (c) $X_2Y_5$}};
\end{tikzpicture}
\caption{Vertex sets $X_2$, $Y_2$ and $X_2Y_5$}
\label{fig:X2Y5}
\end{center}
\end{figure}
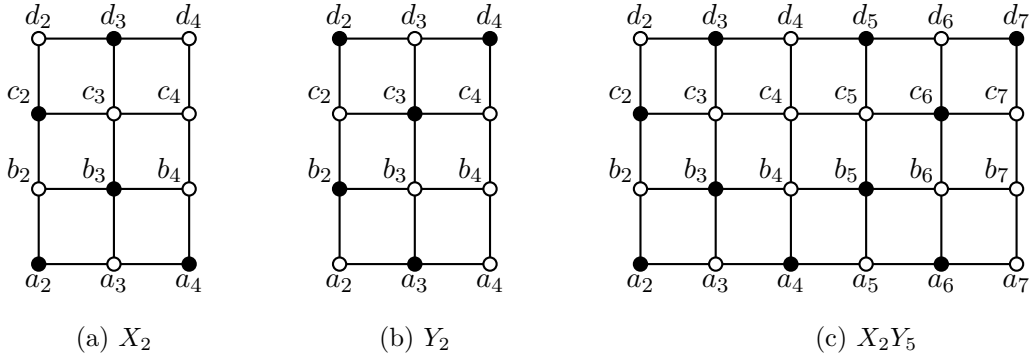

For $i \in \{2,\ldots,m-5\}$, we denote by $X_iY_{i+3}$ the set $X_i \cup Y_{i+3}$, and we denote by $Y_iX_{i+3}$ the set $Y_i \cup X_{i+3}$. The sets $X_2$, $Y_2$ and $X_2Y_5$, for example, are illustrated in Figure~\ref{fig:X2Y5}. We note that all vertices in $V_{i+1}$, $V_{i+2}$, $V_{i+3}$ and $V_{i+4}$ are infected by the set $X_iY_{i+3}$ (respectively, by the set $Y_iX_{i+3}$) in the $4\times 6$ grid induced by the sets $V_i \cup V_{i+1} \cup \cdots \cup V_{i+5}$. For notational simplicity, if the subscripts are clear from the context, we simply write $X$ and $Y$ rather than $X_i$ and $Y_i$, respectively, and we write $XY$ and $YX$ rather than $X_iY_{i+3}$ and $Y_iX_{i+3}$, respectively. We also extend our notation to include multiple copies of $X$ and $Y$. For example, we denote by $X_iY_{i+3}X_{i+6}$ the set $X_i \cup Y_{i+3} \cup X_{i+6}$ and simply denote this by the sequence $XYX$. Using the sequence of sets $XYXY \ldots$, we obtain grids of size $4 \times 3k$ where every internal column from $2$ to $3k-1$ becomes infected. We now construct a percolating set $S$ as follows.

\newpage
\begin{claim}
\label{c:claim15}
If $m \equiv 2 \, (\modo \, 3)$, then $m(G_m,3) \le \lfloor \frac{5(m+1)}{3} \rfloor + 2$.
\end{claim}
\proof Suppose that $m \equiv 2 \, (\modo \, 3)$. Thus, $m = 3k + 2$ for some $k \ge 1$. Let $S$ consist of vertices in $(V_1 \setminus \{c_1\}) \cup V_m$ and from the set $V(G) \setminus (V_1 \cup V_m)$, we add to $S$ the vertices given by the alternating sequence, $XYXY \ldots X$ or $XYXY \ldots Y$, of sets $X$ and $Y$ starting with the set $X$. From our earlier observations, we infer that all vertices in $V_3 \cup V_{m-2}$ become infected. Moreover, the vertices of $S$ infect all vertices in $V_1 \cup V_2$ and infect all vertices in $V_{m-1} \cup V_m$, implying that the set $S$ is a $3$-percolating set. Thus, $m(G_m,3) \le |S| = 7 + 5k = 7 + 5 \times \frac{m-2}{3} = \frac{5}{3}(m+1) + 2$.~\smallqed

\begin{claim}
\label{c:claim16}
If $m \equiv 1 \, (\modo \, 3)$, then $m(G_m,3) \le \lfloor \frac{5(m+1)}{3} \rfloor + 2$.
\end{claim}
\proof Suppose that $m \equiv 1 \, (\modo \, 3)$. Thus, $m = 3k + 1$ for some $k \ge 1$. We now construct the set $S$ as follows. Let $S \cap V_1 = \{a_1,b_1,d_1\}$.

If $m \equiv 1 \, (\modo \, 6)$, then we let $S \cap (V_2 \cup V_3 \cup \cdots \cup V_{m-3})$ consist of the alternating sequence $XYXY \ldots X$ that starts and ends with the set $X$, and we let $S \cap (V_{m-2} \cup V_{m-1} \cup V_{m}) = \{b_{m-2},d_{m-2},a_{m-1},c_{m-1},a_m,b_m,d_m\}$.

If $m \equiv 4 \, (\modo \, 6)$, then we let $S \cap (V_2 \cup V_3 \cup \cdots \cup V_{m-3})$ consist of the alternating sequence $XYXY \ldots Y$ that starts with the set $X$ and ends with the set $Y$, and we let $S \cap (V_{m-2} \cup V_{m-1} \cup V_{m}) = \{a_{m-2},c_{m-2},b_{m-1},d_{m-1},a_m,b_m,d_m\}$.

In both cases, the resulting set $S$ is a $3$-percolating set of $G_m$. Thus, $m(G_m,3) \le |S| =  10 + 5(k-1) = 10 + 5 \times \frac{m-4}{3} = \frac{5}{3}(m+1) + \frac{5}{3}$.~\smallqed

\begin{claim}
\label{c:claim17}
If $m \equiv 0 \, (\modo \, 3)$, then $m(G_m,3) \le \lfloor \frac{5(m+1)}{3} \rfloor + 2$.
\end{claim}
\proof Suppose that $m \equiv 0 \, (\modo \, 3)$. Thus, $m = 3k$ for some $k \ge 2$. We now construct the set $S$ as follows. Let $S \cap V_1 = \{a_1,b_1,d_1\}$.

If $m \equiv 0 \, (\modo \, 6)$, then we let $S \cap (V_2 \cup V_3 \cup \cdots \cup V_{m-2})$ consist of the alternating sequence $XYXY \ldots X$ that starts and ends with the set $X$, and we let $S \cap (V_{m-1} \cup V_{m}) = \{b_{m-1},d_{m-1},a_m,c_m,d_m\}$.

If $m \equiv 3 \, (\modo \, 6)$, then we let $S \cap (V_2 \cup V_3 \cup \cdots \cup V_{m-2})$ consist of the alternating sequence $XYXY \ldots Y$ that starts with the set $X$ and ends with the set $Y$, and we let $S \cap (V_{m-1} \cup V_{m}) = \{a_{m-1},c_{m-1},a_m,b_m,d_m\}$.

In both cases, the resulting set $S$ is a $3$-percolating set of $G_m$. Thus, $m(G_m,3) \le 8 + 5(k-1) = 8 + 5 \times \frac{m-3}{3} + \frac{4}{3}$.~\smallqed

\medskip
By Claims~\ref{c:claim15},~\ref{c:claim16}, and~\ref{c:claim17}, we have $m(G,3) \le \lfloor \frac{5(m+1)}{3} \rfloor + 2$. This completes the proof of Theorem~\ref{thm:4xm-upper}.~\QED

\medskip
Theorem~\ref{thm:4xm} follows as an immediate consequence of Theorems~\ref{thm:4xm-lower} and~\ref{thm:4xm-upper}.

\section{Open problems}

As shown in Theorem~\ref{thm:4xm}, for $m \ge 4$ if $G = P_4 \cp P_m$, then $m(G,3)$ takes on one of two possible values, namely $m(G,3) = \lfloor \frac{5(m+1)}{3} \rfloor + 1$ or $m(G,3) = \lfloor \frac{5(m+1)}{3} \rfloor + 2$. It would be interesting to determine the exact value of $m(G,3)$ in this case for all $m \ge 4$. More generally, it would be interesting to determine the exact value of $m(G,3)$ when $G = P_n \cp P_m$ for all $n \ge m \ge 6$.



\paragraph{Data Availability} Data sharing is not applicable to this article as no new data were created or analyzed in this study.







%
%

\end{document}